\documentclass[journal]{IEEEtran}
\usepackage{latexsym,amsmath,amsfonts,amssymb,mathrsfs,bm,multirow}
\usepackage[ruled,noline]{algorithm2e}
\usepackage{graphicx,epstopdf}
\usepackage{subfigure}
\usepackage{url}
\usepackage{siunitx}
\usepackage{cite}
\usepackage{array}
\usepackage{lineno}
\renewcommand\thesection{\Roman{section}}
\renewcommand\thesubsection{\emph{\Alph{subsection}}}

\newtheorem{thm}{Theorem}
\newtheorem{lem}{Lemma}

\newtheorem{ass}{Assumption}

\ifCLASSINFOpdf
\else
\fi
\allowdisplaybreaks \allowdisplaybreaks[4]
\begin{document}
\title{A variable metric mini-batch proximal stochastic recursive gradient algorithm with diagonal Barzilai-Borwein stepsize}

%\author{Michael~Shell,~\IEEEmembership{Member,~IEEE,}
%        John~Doe,~\IEEEmembership{Fellow,~OSA,}
%        and~Jane~Doe,~\IEEEmembership{Life~Fellow,~IEEE}% <-this % stops a space
%\thanks{M. Shell was with the Department
%of Electrical and Computer Engineering, Georgia Institute of Technology, Atlanta,

%\author{Tengteng~Yu, Xin-Wei~Liu, Yu-Hong~Dai and~Jie~Sun}% <-this % stops a space

\author{Tengteng~Yu, Xin-Wei~Liu, Yu-Hong~Dai and~Jie~Sun% <-this % stops a space

\thanks{Tengteng Yu is with the School of Artificial Intelligence, Hebei University of Technology, Tianjin 300401, China (e-mail: yuteng206@163.com).}
\thanks{Xin-Wei Liu is with the Institute of Mathematics, Hebei University of Technology, Tianjin 300401, China (e-mail: mathlxw@hebut.edu.cn).}
\thanks{Yu-Hong Dai is with the State Key Laboratory of Scientific and Engineering Computing (LSEC), Academy of Mathematics and Systems Science, Chinese Academy of Sciences, Beijing 100190, China, and also with the School of Mathematical Sciences, University of Chinese Academy of Sciences, Beijing 100049, China (e-mail: dyh@lsec.cc.ac.cn).}
\thanks{Jie Sun is s with the Institute of Mathematics, Hebei University of Technology, Tianjin 300401, China, and also with the School of Business, National University of Singapore, Singapore 119245 (e-mail:  jie.sun@nus.edu.sg).}
}

\maketitle

\begin{abstract}
Variable metric proximal gradient methods with different metric selections have been widely used in composite optimization. Combining the Barzilai-Borwein (BB) method with a diagonal selection strategy for the metric, the diagonal BB stepsize can keep low per-step computation cost as the scalar BB stepsize and better capture the local geometry of the problem. In this paper, we propose a variable metric mini-batch proximal stochastic recursive gradient algorithm VM-mSRGBB, which updates the metric using a new diagonal BB stepsize. The linear convergence of VM-mSRGBB is established for strongly convex, non-strongly convex and convex functions. Numerical experiments on standard data sets show that VM-mSRGBB is better than or comparable to some variance reduced stochastic gradient methods with best-tuned scalar stepsizes or BB stepsizes. Furthermore, the performance of VM-mSRGBB is superior to some advanced mini-batch proximal stochastic gradient methods.
\end{abstract}


% \fi
\IEEEpeerreviewmaketitle

%\section{Introduction}
%% The very first letter is a 2 line initial drop letter followed
%% by the rest of the first word in caps.
%%
%% form to use if the first word consists of a single letter:
%% \IEEEPARstart{A}{demo} file is ....
%%
%% form to use if you need the single drop letter followed by
%% normal text (unknown if ever used by the IEEE):
%% \IEEEPARstart{A}{}demo file is ....
%%
%% Some journals put the first two words in caps:
%% \IEEEPARstart{T}{his demo} file is ....
%%
%% Here we have the typical use of a "T" for an initial drop letter
%% and "HIS" in caps to complete the first word.

\section{Introduction}\label{sin}
\IEEEPARstart{W}{e} consider the following problem of minimizing a composition of two convex functions:
\begin{equation}\label{eqvm1.1}
  \min_{w\in \mathbb{R}^d}  P(w)= F(w)+ R(w),
\end{equation}
where $F(w)=\frac{1}{n}\sum_{i=1}^{n}f_i(w)$, each component function $f_i(w):\mathbb{R}^d\rightarrow \mathbb{R}$, $i=1,2,\ldots,n$ is convex and smooth, $n$ is the sample size, and $R(w): \mathbb{R}^d\rightarrow \mathbb{R}\cup\{+\infty\}$ is a relatively simple proper convex function and can be non-differentiable.  The term $R(w)$ is sometimes referred to as a regularization. In this paper, we are especially interested in the case where $n$ is extremely large, and the proximal operator of $R(w)$ can be computed efficiently.

The formulation \eqref{eqvm1.1} appears across a broad range of applications in machine learning \cite{sra2012Optimization,Shalev2014Understanding,bottou2018optimization}, statistics \cite{hastie2009elements}, matrix completion \cite{recht2013parallel}, neural networks \cite{zhang2015deep,Goodfellow2016Deep,li2017preconditioned}, etc. One popular instance is the regularized empirical risk minimization (ERM) \cite{bottou2018optimization,hastie2009elements,jin2018stochastic}, which involves a collection of training examples $\{(a_i,b_i)\}_{i=1}^n$, where $a_i\in\mathbb{R}^d$ is a feature vector and $b_i\in\mathbb{R}$ is the desired response. With the component functions $f_i(w) = \frac{1}{2}(b_i-a_i^Tw)$, Lasso, ridge regression and elastic net employ the regularization terms $R(w) = \lambda_1\|w\|_1$, $R(w) = \frac{\lambda_2}{2}\|w\|_2$ and $R(w) = \lambda_1\|w\|_1+\frac{\lambda_2}{2}\|w\|_2$, respectively, where $\lambda_1$ and $\lambda_2$ are nonnegative regularization parameters. When considering binary classification problems, one frequently used component function is the logistic loss $f_i(w) = \log(1+\exp(-b_ia_i^Tw))$ and $R(w)$ can be any of the above regularization terms.

One of the most popular methods for solving optimization problems in composite form \eqref{eqvm1.1} is the proximal gradient descent (Prox-GD), which has attracted many researchers in improving computation costs, establishing theoretical convergence results under mild conditions, and designing practical rules for stepsize selections \cite{combettes2005signal,zhou2006gradient,beck2009fast}. Some accelerated Prox-GD variants have also been proposed, see for example \cite{nesterov2013gradient,parikh2014proximal,drusvyatskiy2018optimal,bubeck2015geometric}. However, problem \eqref{eqvm1.1} with a large sum of $n$ component functions becomes challenging for Prox-GD since it requires computing the exact full gradient. Motivated by the seminal work of Robbins and Monro \cite{robbins1951stochastic}, a proximal stochastic gradient descent (Prox-SGD) method has been developed, which chooses $i_k\in\{1,2,\ldots,n\}$ uniformly at random and takes the update
\begin{equation}\label{eqvm1.2}
w_{k+1} = \arg\min_{w\in \mathbb{R}^d} \{\nabla f_{i_k}(w_k)^Tw + \frac{1}{2\eta_k}\|w-w_k\|_2^2 + R(w)\},
\end{equation}
where $\nabla f_{i_k}(w_k)$ is the gradient of the $i_k$-th component function $f_{i_k}$ at $w_k$ and $\eta_k>0$ is the stepsize (a.k.a. learning rate). Let us define the scaled proximal operator of $R$ relative to the metric $A$ \cite{park2019variable} by
\begin{equation}\label{eqvmpm}
  \textrm{prox}_{R}^A(w) = \arg\min_{y\in \mathbb{R}^d}\{\frac{1}{2}\|y-w\|_A^2 + R(y)\},
\end{equation}
where $A\in \mathbb{R}_{++}^{d\times d}$ is a positive definite matrix and $\|z\|_A = \sqrt{z^TAz}$ is the  norm induced by $A$ (or $A$-norm), then the update rule of Prox-SGD can be described more compactly as
\begin{equation}\label{eqvm1.3}
 w_{k+1} = \textrm{prox}_{R}^{\eta_k^{-1} I}\big(w_k - \eta_k \nabla f_{i_k}(w_k)\big),
\end{equation}
where $I\in \mathbb{R}^{d\times d}$ is the identity matrix. When $R(w)$ is a constant function, the update rule in \eqref{eqvm1.3} becomes the standard SGD method.

Prox-SGD has the great advantage of  tremendous per-iteration saving since it evaluates the gradient of a single component function rather than the full gradient. Due to the large variance of the stochastic gradient introduced by random sampling, Prox-SGD only enjoys a sublinear convergence rate for strongly convex functions as opposed to a linear convergence rate of Prox-GD. Starting from several prevalent variance reduced stochastic gradient methods such as SAG \cite{roux2012a,schmidt2017minimizing}, SVRG \cite{johnson2013accelerating}, SAGA \cite{defazio2014saga:}, S2GD \cite{konevcny2013semi}, SARAH \cite{nguyen2017sarah} and SPIDER \cite{fang2018spider:}, recent works consider to incorporate variance reduction techniques to improve the convergence rate of Prox-SGD. In \cite{xiao2014proximal}, Xiao and Zhang proposed a proximal variant of SVRG, called Prox-SVRG, and proved its linear convergence rate for strongly convex problems. By combining mini-batch scheme with S2GD, Kone\v{c}n\'{y} et al. \cite{konevcny2015mini} developed the mS2GD method that achieves better theoretical complexity and practical performance than Prox-SVRG. A proximal version of SARAH can be found in \cite{pham2020proxsarah}.

Since the stepsize has an important influence on the performances of stochastic gradient methods, many researchers are devoted to designing more efficient scheme of stepsizes. For classical SGD, one frequently employed stepsize strategy in practical computation is
\begin{equation*}
  \sum_{k=1}^{\infty} \eta_k = \infty \quad \textrm{and } \quad \sum_{k=1}^{\infty} \eta_k^2 < \infty.
\end{equation*}
However, such a choice often yields sublinear convergence of SGD, see \cite{bottou2018optimization} for example. In recent years, using the Barzilai-Borwein (BB) method \cite{barzilai1988two} to automatically calculate stepsizes for SGD and its variants has attracted more and more attention. One great advantage of BB stepsize is that it is able to capture hidden second order information and is insensitive to the choice of initial stepsizes, which makes it very promising in practice. See \cite{barzilai1988two,dhl2018,fletcher2005barzilai} and references therein for more details about BB-like methods. One pioneer work in this line is due to Tan et al. \cite{tan2016barzilai}, who proposed to incorporate the BB stepsize with SGD and SVRG, and got the SGD-BB and SVRG-BB methods. By combining SARAH with the BB method and importance sampling strategy, Liu et al. \cite{liu2020a} suggested the SARAH-I-BB method. To solve problem \eqref{eqvm1.1}, Yu et al. \cite{yu2019trust} developed a mini-batch proximal stochastic recursive gradient algorithm that incorporates the trust-region scheme and BB stepsize.

Recently, Park et al. \cite{park2019variable} proposed a variable metric proximal gradient method, called VM-PG, for minimizing composite functions, which uses an adaptive metric selection strategy called the diagonal BB stepsize. As pointed out in \cite{park2019variable}, the diagonal BB stepsize can better capture the local geometry of the problem and keep per-step computation cost similar to the scalar BB stepsize. However, VM-PG is designed in the deterministic form and cannot be directly applied to large-scale machine learning problems \cite{bonnans1995a,bonettini2016a,salzo2017the}.

In this paper, motivated by VM-PG and the success of SGD and its variants in solving problem \eqref{eqvm1.1}, we propose a mini-batch proximal stochastic recursive gradient method, named VM-mSRGBB. The proposed VM-mSRGBB method updates the metric by a new diagonal BB stepsize, which is the closed-form solution of a constrained optimization. In each iteration, the computational cost on gradients of our VM-mSRGBB method is the same as that of SVRG and SARAH.
%to use the diagonal BB stepsize to adjust the metric for the mini-batch proximal stochastic recursive gradient method, which yields a new method VM-mSRGBB. The per-iteration cost of VM-mSRGBB is almost the same as SARAH.
We present the convergence analysis of VM-mSRGBB under different conditions, which shows that it converges linearly for strongly convex, non-strongly convex and convex functions. Numerical results for solving regularized logistic regression problems on standard data sets show that the performance of VM-mSRGBB is better than or comparable to Prox-SVRG with best-tuned stepsizes and the proximal variant of SVRG-BB with different initial stepsizes. Further comparisons between VM-mSRGBB and some advanced mini-batch proximal stochastic gradient methods demonstrate the efficiency of VM-mSRGBB.

The rest of this paper is organized as follows. In Section \ref{vmmsrgbb} we propose our VM-mSRGBB method. In Section \ref{converge} we prove that VM-mSRGBB always enjoy a linear convergence rate under different conditions. Numerical experiments are then reported in Section \ref{experiment}. Finally, we draw some conclusions in Section \ref{conclusion}.

\section{The VM-mSRGBB method}\label{vmmsrgbb}

Our VM-mSRGBB method is motivated by the VM-PG method for solving composite problems in the deterministic setting, which uses a variable metric rather than a scalar matrix to estimate the second-order information of $F(w)$ and provides better approximation of the local Hessian at each step. A formal description of VM-mSRGBB is given in Algorithm \ref{conVM}.
\begin{algorithm}
\DontPrintSemicolon
\caption{VM-mSRGBB($\tilde{w}^0,m,b,U_0$)}\label{conVM}

\textbf{Input:} update frequency $m$ (max \# of stochastic steps per outer loop), initial point $\tilde{w}^0 \in \mathbb{R}^d$, initial matrix $U_0 = \eta_0I$, mini-batch size $b\in \{1,2,\ldots,n\}$;

\For{$k=0,1,\ldots,K-1$}{

$w_1^{k} = w_0^{k} = \tilde{w}^k$;

$v_0^{k} = \nabla F(w_0^k)$;

Probability $Q = \{q_1,q_2,\ldots,q_n\}$ on $\{1,2,\ldots,n\}$;

Choose $t_k\in\{1,2,\ldots,m\}$ uniformly at random;

\For{$t=1,\ldots,t_k$}{

Choose mini-batch $I_t\subseteq \{1,2,\ldots, n\}$ of size $b$, where each $i\in I_t$ is chosen from $\{1,2,\ldots, n\}$ randomly according to $Q$;

\begin{equation}\label{eqvmvt1}
  v_t^k =  \frac{1}{b}\sum_{i\in I_t}\big[(\nabla f_i(w_t^k) - \nabla f_i(w_{t-1}^k))/(q_in)\big] + v_{t-1}^k;
\end{equation}

$w_{t+1}^k = \textrm{prox}_{R}^{U_k^{-1}}(w_t^k- U_kv_t^{k})$;

}

$\tilde{w}^{k+1} =  w_{t_k+1}^k$;

Compute $U_{k}$ from \eqref{eqvmstep};

}

\textbf{Output:} Iterate $w_a$ chosen uniformly at random from $\{\{w_t^k\}_{t=1}^{t_k}\}_{k=0}^{K-1}$;
\end{algorithm}

Before presenting the selection of the metric $U_k$, we would like to mention that $v_t^k$ is a biased estimate of the full gradient $\nabla F(w_t^k)$, which is the same as SARAH \cite{nguyen2017sarah} but different from SGD and SVRG types of methods \cite{johnson2013accelerating,xiao2014proximal}. In fact, it is easy to see that the conditional expectation of $v_t^k$ given $\mathcal{F}_t$ is
\begin{eqnarray*}\label{eqvmtotal1}
\mathbb{E}[v_t^k|\mathcal{F}_t] &=&  \sum_{i=1}^n \frac{\nabla f_i(w_t^k) - \nabla f_i(w_{t-1}^k)}{q_in}\cdot q_i + v_{t-1}^k \nonumber\\
&=& \nabla F(w_t^k) - \nabla F(w_{t-1}^k) + v_{t-1}^k,
\end{eqnarray*}
where $\mathcal{F}_t = \sigma(w_0^k,I_1,I_2,\ldots,I_{t-1})$ is the $\sigma$-algebra generated by $w_0^k,I_1,I_2,\ldots,I_{t-1}$ and $\mathcal{F}_0 = \mathcal{F}_1 = \sigma(w_0^k)$.
%
%we can easily show $v_t^k$ is a biased estimate by taking expectation with respect to $I_t$. Let us denote $\mathcal{F}_t = \sigma(w_0^k,I_1,I_2,\ldots,I_{t-1})$ by the $\sigma$-algebra generated by $w_0^k,I_1,I_2,\ldots,I_{t-1}$ and $\mathcal{F}_0 = \mathcal{F}_1 = \sigma(w_0^k)$. For any $t\geq 1$, conditioned on $\mathcal{F}_t$, we obtain
%\begin{eqnarray*}\label{eqvmtotal1}
%% \nonumber to remove numbering (before each equation)
%\mathbb{E}[v_t^k|\mathcal{F}_t] &=&  \sum_{i=1}^n \frac{\nabla f_i(w_t^k) - \nabla f_i(w_{t-1}^k)}{q_in}\cdot q_i + v_{t-1}^k \nonumber\\
%&=& \nabla F(w_t^k) - \nabla F(w_{t-1}^k) + v_{t-1}^k.
%\end{eqnarray*}
As will be seen in Theorems \ref{vmthmdes} and \ref{vmthm2.1}, the simple recursive framework for updating $v_t^k$ yields a non-increasing property and a linear convergence of the inner loop of our VM-mSRGBB method, which does not hold for Prox-SVRG and mS2GD.

When taking total expectation and employing the fact $v_0^k = \nabla F(w_0^k)$, it follows that $\mathbb{E}[v_1^k] = \mathbb{E}[\nabla F(w_1^k)] - \mathbb{E}[\nabla F(w_0^k)] + \mathbb{E}[v_0^k] = \mathbb{E}[\nabla F(w_1^k)]$. By induction, we obtain
\begin{equation}\label{eqvmtotal2}
\mathbb{E}[v_t^k] = \mathbb{E}[\nabla F(w_t^k)].
\end{equation}

Notice that, when $U_k = \alpha_kI$ with $\alpha_k$ being a scalar stepsize,  Algorithm \ref{conVM} is a proximal version of SARAH \cite{nguyen2017sarah}. And it transforms to the stochastic proximal quasi-Newton method for $U_k\approx (\nabla^2 F(w_t^k))^{-1}$ \cite{wang2017inexact,wang2019stochastic}. However, a scalar stepsize cannot capture the inverse Hessian well and the inverse Hessian may be expensive to calculate.
%Moreover, on the one hand, the diagonal structure used in VM-mSRGBB can better satisfy the secant condition than a single scalar element used in BB method. On the other hand, the existence of $\omega$ makes it robust to the degenerate case where $s_k^Ty_k \approx 0$. These makes VM-mSRGBB more stable than that used the scalar BB stepsizes.
Motivated by \cite{park2019variable}, we suggest a diagonal metric $U_k$ computed as follows
\begin{eqnarray}\label{eqvmvm}
&\min_{u\in \mathbb{R}^d}& \|s_k-Uy_k\|_2^2 + \omega \|U-U_{k-1}\|_F^2 \\
&\textrm{s.t.}& \alpha_k^2I \preceq U \preceq \alpha_k^1I, \nonumber\\
&& U = \textrm{Diag}(u), \nonumber
\end{eqnarray}
where $s_k = \tilde{w}^k - \tilde{w}^{k-1}$, $y_k = \nabla F(\tilde{w}^k) - \nabla F(\tilde{w}^{k-1})$, $\|\cdot\|_F$ is the Frobenius norm and $0< \alpha_k^2\leq \alpha_k^1$ are two stepsizes given by users. Clearly, the solution $U_k$ of \eqref{eqvmvm} satisfies the secant equation $s_k=U_ky_k$ in the sense of least squares and is close to the previous metric $U_{k-1}$ where the closeness is controlled by the hyperparameter $\omega> 0$. So, $U_k$ can capture the geometry of the inverse Hessian of $F(w)$, which is different from the one in \cite{park2019variable}.

%Notice that the hyperparameter $\omega> 0$ controls the trade-off between satisfying the secant equation and being consistent with the previous metric $U_{k-1}$.

For $U_k = \textrm{Diag}(u_k)\in \mathbb{R}^{d\times d}$ with $u_k= [u_k^1,u_k^2,\ldots,u_k^d]\in \mathbb{R}^{d}$, problem \eqref{eqvmvm} has a closed-form solution given by
\begin{equation}\label{eqvmstep}
u_k^i = \left\{
\begin{array}{ll}
\alpha_k^2, & \hbox{$\frac{s_k^iy_k^i + \omega u_{k-1}^i}{(y_k^i)^2 + \omega} < \alpha_k^2$;} \\
\alpha_k^1, & \hbox{$\frac{s_k^iy_k^i + \omega u_{k-1}^i}{(y_k^i)^2 + \omega} > \alpha_k^1$;} \\
\frac{s_k^iy_k^i + \omega u_{k-1}^i}{(y_k^i)^2 + \omega}, & \hbox{otherwise.}
\end{array}
\right.
\end{equation}
where $s_k^i$ and $y_k^i$ are the $i$-th elements of $s_k$ and $y_k$, respectively.

%In order to save time on choosing the two stepsizes $\alpha_k^1$ and $\alpha_k^2$, we would like to employ BB-like methods to compute stepsizes automatically, which has been successfully used in SGD and its variants, such as SGD-BB and SVRG-BB \cite{tan2016barzilai}.

As mentioned before, the BB stepsize is suitable for SGD and its variants. We would like to employ BB-like stepsizes for $\alpha_k^1$ and $\alpha_k^2$. Since at most $m$ biased gradient estimators are added to $w_0^k$ for getting $w_m^k$ in the inner loop, we employ the following stepsizes
\begin{equation}\label{BB1}
  \alpha_k^{1}= \frac{2}{m}\cdot\frac{\|s_{k}\|_2}{\|y_{k}\|_2}
\end{equation}
and
\begin{equation}\label{BB2}
\alpha_k^{2}= \frac{1}{m}\cdot\frac{s_{k}^Ty_{k}}{\|y_{k}\|_2^2}.
\end{equation}
Here, $\alpha_k^{1}$ is a variant of the BB-like stepsize $\alpha_k^D=\frac{\|s_{k}\|_2}{\|y_{k}\|_2}$ proposed in \cite{dai2015positive} and $\alpha_k^{2}$ is a variant of the original BB stepsize $\alpha_k^{BB}=\frac{s_{k}^Ty_{k}}{\|y_{k}\|_2^2}$ in \cite{barzilai1988two}.
%That is, $\alpha_k^{1}$ is $2/m$ of the positive stepsize $\alpha_k^D=\frac{\|s_{k}\|_2}{\|y_{k}\|_2}$ and $\alpha_k^{2}$ equals to the BB stepsize  $\alpha_k^{BB}=\frac{s_{k}^Ty_{k}}{\|y_{k}\|_2^2}$ divided by $m$.
Notice that by the Cauchy-Schwarz inequality $\alpha_k^D\geq\alpha_k^{BB}$ always holds. Moreover, $\alpha_k^D$ can be seen as an approximation of $1/L$ with $L$ being the Lipschitz constant of $\nabla F$, see \cite{dai2015positive}.

%In \cite{nesterov1998introductory}, the authors show that when the stepsize is chosen from $(0,2/L)$, the gradient method converges linearly for strongly convex functions.

\section{Convergence analysis}\label{converge}
In order to establish convergence of VM-mSRGBB in different cases, we make the following two blanket assumptions.
\begin{ass}\label{vmnonsmooth}
The regularization function $R(w):\mathbb{R}^d\rightarrow \mathbb{R}\cup\{+\infty\}$ is a lower semi-continuous and convex function. However, it can be non-differentiable. Its effective domain, $\textrm{dom}(R) = \{w\in \mathbb{R}^d| R(w)< + \infty\}$, is closed.
\end{ass}

\begin{ass}\label{vmsmooth}
Each component function $f_i(w):\mathbb{R}^d \rightarrow \mathbb{R}$ is convex and $L_i$-smooth, that is, there exists $L_i>0$ such that
\begin{equation*}
  \|\nabla f_i(w) - \nabla f_i(w')\|_2 \leq L_i\|w-w'\|_2, \quad \forall w, w'\in \textrm{dom}(R).
\end{equation*}
\end{ass}

Assumption \ref{vmsmooth} implies that $F(w)$ is also $L$-smooth with $L \leq \frac{1}{n}\sum_{i=1}^n L_i$. For simplicity, we denote $L_{\Omega}$ as
\begin{equation*}
  L_{\Omega} = \max_{i=1,2,\ldots,n} \frac{L_i}{nq_i},
\end{equation*}
then $L_{\Omega} \geq \frac{1}{n}\sum_{i=1}^n L_i \geq L$. It is not difficult to obtain the following  result from Assumption \ref{vmsmooth}.
\begin{lem}\label{vmlemAs1}(Theorem 2.1.5 \cite{nesterov1998introductory})
Suppose that $f_i$ is convex and $L_i$-smooth. Then, for any $w, w'\in \mathbb{R}^d$,
\begin{equation*}\label{eqvmAs1}
  (\nabla f_i(w) - \nabla f_i(w'))^T(w-w') \geq\frac{1}{L_i}\|\nabla f_i(w) - \nabla f_i(w')\|_2^2.
\end{equation*}
\end{lem}

Now we generalize some basic properties of proximal mapping to scaled proximal operator. Although they are direct extensions, we have not find the same results in literature.
\begin{lem}\label{vmlem2}
	Let $R(w)$ be a proper closed and convex function on $\mathbb{R}^d$. Then $\textrm{prox}_{R}^{A^{-1}}(w)$ is a singleton for any $w\in \textrm{dom}(R)$ and any symmetric positive definite matrix $A\in \mathbb{S}_{++}^{d\times d}$. Furthermore, the following statements are equivalent:
\begin{itemize}
	\item [(i)] $\mathbf{u}=\textrm{prox}_{R}^{A^{-1}}(w)$.
	\item [(ii)] $A^{-1}(w-\mathbf{u}) \in \partial R(\mathbf{u})$, where $\partial R$ is the subdifferential of $R$.
\end{itemize}
\end{lem}
\begin{IEEEproof}
	The uniqueness of $\textrm{prox}_{R}^{A^{-1}}(w)$ can be proved in a similar way as Theorem 6.3 of \cite{Beck2017First} by noting that $A$ is symmetric positive definite. For the latter part, one can employ the techniques used in the proof of Theorem 6.39 in \cite{Beck2017First}. We omit the details here.
\end{IEEEproof}

\begin{lem}\label{vmlem6}
Let $R(w)$ be a proper closed and convex function on $\mathbb{R}^d$. Then, for any $w,w'\in \textrm{dom}(R)$ and any $A\in \mathbb{S}_{++}^{d\times d}$, the following inequality holds:
\begin{equation*}
  \|\textrm{prox}_{R}^{A^{-1}}(w)-\textrm{prox}_{R}^{A^{-1}}(w')\|_{A^{-1}}^2 \leq \|w-w'\|_{A^{-1}}^2.
\end{equation*}
\end{lem}

\begin{IEEEproof}
See Appendix \ref{vmapp4}.
\end{IEEEproof}

The following theorem shows that our proximal stochastic recursive step $w_{t+1}^k - w_t^k$ decreases in expectation for convex functions.
\begin{thm}\label{vmthmdes}
Suppose that Assumptions \ref{vmnonsmooth} and \ref{vmsmooth} hold. Consider $v_t^k$ defined by \eqref{eqvmvt1} in VM-mSRGBB (Algorithm \ref{conVM}) with $0 \prec U_k \preceq 1/L_{\Omega}I$. Then, in the $k$-th outer loop, for any $t> 1$, we have
\begin{equation*}
  \mathbb{E}\big[\|w_{t+1}^k - w_t^k\|_{U_k^{-1}}^2\big] \leq   \mathbb{E}\big[\|w_t^k - w_{t-1}^k\|_{U_k^{-1}}^2\big],
\end{equation*}
where the expectation is taken with respect to all the variables generated in the $k$-th outer loop.
\end{thm}

\begin{IEEEproof}
We take expectation on $\|w_{t+1}^k - w_t^k\|_{U_k^{-1}}^2$ with respect to all the variables generated in the $k$-th outer loop and obtain
\begin{align}\label{eqvmc1.1}
  & \mathbb{E}\big[\|w_{t+1}^k - w_t^k\|_{U_k^{-1}}^2\big] \nonumber\\
  =& \mathbb{E}\big[\|\textrm{prox}_{R}^{U_k^{-1}}(w_t^k - U_k v_t^k) - \textrm{prox}_{R}^{U_k^{-1}}(w_{t-1}^k - U_k v_{t-1}^k)\|_{U_k^{-1}}^2 \big] \nonumber\\
  \leq& \mathbb{E}\big[\|w_t^k - w_{t-1}^k -U_k(v_t^k - v_{t-1}^k)\|_{U_k^{-1}}^2\big] \nonumber\\
  =& \mathbb{E}\big[\|w_t^k - w_{t-1}^k\|_{U_k^{-1}}^2 + \|v_t^k - v_{t-1}^k\|_{U_k}^2\big] \nonumber\\
  & - \mathbb{E}\big[2(w_t^k - w_{t-1}^k)^T(v_t^k - v_{t-1}^k)\big] \nonumber\\
  =&\mathbb{E}\big[\|w_t^k - w_{t-1}^k\|_{U_k^{-1}}^2\big] + \mathbb{E}\big[\|v_t^k - v_{t-1}^k\|_{U_k}^2\big] \nonumber\\
  & - 2\mathbb{E}\big[(w_t^k - w_{t-1}^k)^T(\frac{1}{b}\sum_{i\in I_t}\frac{\nabla f_{i}(w_t^k)-\nabla f_{i}(w_{t-1}^k)}{q_{i}n})\big] \nonumber\\
  \leq& \mathbb{E}\big[\|w_t^k - w_{t-1}^k\|_{U_k^{-1}}^2\big] + \mathbb{E}\big[\|v_t^k - v_{t-1}^k\|_{U_k}^2\big] \nonumber\\
  & - 2\mathbb{E}\big[\frac{1}{b}\sum_{i\in I_t}\frac{\|\nabla f_{i}(w_t^k)-\nabla f_{i}(w_{t-1}^k)\|_2^2}{q_{i}nL_i}\big] \nonumber\\
  \leq & \mathbb{E}\big[\|w_t^k - w_{t-1}^k\|_{U_k^{-1}}^2\big]  + \mathbb{E}\big[\|v_t^k - v_{t-1}^k\|_{U_k}^2\big] \nonumber\\
  & - \frac{2}{L_{\Omega}}\mathbb{E}\big[\frac{1}{b}\sum_{i\in I_t}\|\frac{\nabla f_{i}(w_t^k)-\nabla f_{i}(w_{t-1}^k)}{q_{i}n}\|_2^2\big] \\
  \leq &\mathbb{E}\big[\|w_t^k - w_{t-1}^k\|_{U_k^{-1}}^2\big]  + \frac{1}{L_{\Omega}} \mathbb{E}\big[\|v_t^k - v_{t-1}^k\|_2^2\big] \nonumber\\
  & - \frac{2}{L_{\Omega}}\mathbb{E}\big[\|\frac{1}{b}\sum_{i\in I_t}\frac{\nabla f_{i}(w_t^k)-\nabla f_{i}(w_{t-1}^k)}{q_{i}n}\|_2^2\big] \nonumber\\
   =& \mathbb{E}\big[\|w_t^k - w_{t-1}^k\|_{U_k^{-1}}^2\big] - \frac{1}{L_{\Omega}}\mathbb{E}\big[\|v_t^k - v_{t-1}^k\|_2^2\big] \nonumber\\
  \leq & \mathbb{E}\big[\|w_t^k - w_{t-1}^k\|_{U_k^{-1}}^2\big],\nonumber
\end{align}
where the first inequality follows from Lemma \ref{vmlem6}, and the second inequality uses Lemma \ref{vmlemAs1}. The third inequality holds due to $L_{\Omega} \geq L_i/(nq_i)$ for $i=1,2,\ldots,n$. In the fourth inequality we use the fact that $\mathbb{E}\big[\|z_1+z_2+\ldots+z_r\|_2^2\big] \leq r\mathbb{E}\big[\|z_1\|_2^2 + \|z_2\|_2^2 + \ldots + \|z_r\|_2^2\big]$ with $z_j$ being random variables for $j\in\{1,2,\ldots,r\}$ and  $0\prec U_k \preceq 1/L_{\Omega}I$. The last equality holds by the definition of $v_t^k$.
\end{IEEEproof}

Let $\mathcal{W}_*$ be the set of optimal solutions of problem \eqref{eqvm1.1} and $w_*\in\mathcal{W}_*$. From Theorem 2 in \cite{yu2019trust}, an upper bound on the variance of $v_t^k$ can be given as follows.
\begin{lem}\label{vmlem1}
Suppose that Assumptions \ref{vmnonsmooth} and \ref{vmsmooth} hold, and choose $b\in\{1,2,\ldots,n\}$. Consider $v_t^k$ as defined in \eqref{eqvmvt1}. Then, for $t = 1, 2,\ldots,m$, we have
\begin{align*}
  & \mathbb{E}\big[\|v_t^k - \nabla F(w_t^k)\|_2^2\big] \\
  \leq& \frac{4L_{\Omega}}{b}\mathbb{E}\big[P(w_t^k) - P(w_*)
   + P(w_{t-1}^k) - P(w_*)\big],
\end{align*}
where the expectation is taken with respect to all the variables
generated in the $k$-th outer loop.
\end{lem}

To analyze the convergence of multiple outer loops, we define the following generalization of stochastic gradient mapping:
\begin{equation}\label{eqvmspm}
  g_t^k = U_k^{-1}(w_t^k - w_{t+1}^k) = U_k^{-1}\big(w_t^k - \textrm{prox}_{R}^{U_k^{-1}}(w_t^k - U_kv_t^k)\big).
\end{equation}
Then the proximal stochastic gradient step in Algorithm \ref{conVM} can be written as
\begin{equation}\label{eqvmspm2}
  w_{t+1}^k = w_t^k - U_kg_t^k.
\end{equation}

Before establishing the convergence of VM-mSRGBB, we show an upper bound on $P(w)$ by using \eqref{eqvmspm} and \eqref{eqvmspm2} in a similar way to Lemma 3.7 in \cite{xiao2014proximal}. However, we do not require the strong convexity of $F(w)$ and $R(w)$.
\begin{lem}\label{vmlem5}
Suppose that Assumptions \ref{vmnonsmooth} and \ref{vmsmooth} hold, and $0\prec U_k \preceq 1/L_{\Omega} I$. For any $t\geq1$, we have
\begin{equation*}
(w_*-w_t)^T g_t^k + \frac{1}{2}\|g_t^k\|_{U_k}^2 \leq P(w_*) - P(w_{t+1}^k) - (w_*-w_{t+1}^k)^T\delta_t^k,
\end{equation*}
where $\delta_t^k = \nabla F(w_t^k)-v_t^k$.
\end{lem}

\begin{IEEEproof}
See Appendix \ref{vmapp1}.
\end{IEEEproof}

\subsection{VM-mSRGBB for strongly convex functions}
We analyze the linear convergence of VM-mSRGBB in the case where $P(w)$ is strongly convex.
\begin{ass}\label{vmstrong}
The objective function $P(w)$ is $\mu$-strongly convex, that is, there exits $\mu>0$ such that for all $w\in\textrm{dom}(R)$ and $w'\in\mathbb{R}^d $,
\begin{equation*}
  P(w') \geq P(w) +  \xi^T(w-w') + \frac{\mu}{2}\|w-w'\|_2^2, \textrm{  $\forall \xi \in \partial P(w)$}.
\end{equation*}
\end{ass}
Either $F(w)$ or $R(w)$ or both may bring about the strong convexity of $P(w)$. Assumptions \ref{vmnonsmooth}, \ref{vmsmooth} and \ref{vmstrong} are often satisfied by objective functions in machine learning, such as ridge regression and elastic net regularization logistic regression. Moreover, $w_*$ is unique when $P(w)$ is strongly convex.

%Consider $v_t^k$ defined by \eqref{eqvmvt1} in VM-mSRGBB (Algorithm \ref{conVM}).
The following theorem shows that our proximal stochastic recursive step has a linear convergence rate for strongly convex functions.
\begin{thm}\label{vmthm2.1}
Suppose that Assumptions \ref{vmnonsmooth} and \ref{vmsmooth} hold, $F(w)$ is $\mu_F$-strongly convex and $0\prec U_k \preceq 2/L_{Q}I$. Then, in the $k$-th outer loop, for any $t> 1$, we have
\begin{align*}
  & \mathbb{E}\big[\|w_{t+1}^k - w_t^k\|_{U_k^{-1}}^2\big] \\
  \leq&  \big(1- (\mu_F^2 u_k^{\min})(\frac{2}{L_{\Omega}} - u_k^{\max})\big) \mathbb{E}\big[\|w_t^k - w_{t-1}^k\|_{U_k^{-1}}^2\big],
\end{align*}
where $u_k^{\max} = \max_j\{u_k^{j}\}$ , $u_k^{\min} = \min_j\{u_k^{j}\}$ and the expectation is taken with respect to all the variables generated in the $k$-th outer loop.
\end{thm}
\begin{IEEEproof}
The inequality \eqref{eqvmc1.1} in Theorem \ref{vmthmdes} indicates that
\begin{align*}
  & \mathbb{E}\big[\|w_{t+1}^k - w_t^k\|_{U_k^{-1}}^2\big]\\
  \leq & \mathbb{E}\big[\|w_t^k - w_{t-1}^k\|_{U_k^{-1}}^2\big] +(u_k^{\max}-\frac{2}{L_{\Omega}}) \mathbb{E}\big[\|v_t^k - v_{t-1}^k\|_2^2\big]\\
  \leq & \mathbb{E}\big[\|w_t^k - w_{t-1}^k\|_{U_k^{-1}}^2\big]\\
   &+(u_k^{\max}-\frac{2}{L_{\Omega}}) \mathbb{E}\big[\|\nabla F(w_t^k) - \nabla F(w_{t-1}^k)\|_2^2\big]\\
  \leq & \mathbb{E}\big[\|w_t^k - w_{t-1}^k\|_{U_k^{-1}}^2\big] \\
  &+ \mu_F^2(u_k^{\max}-\frac{2}{L_{\Omega}}) \mathbb{E}\big[\|w_t^k - w_{t-1}^k\|_2^2\big]\\
  \leq & \mathbb{E}\big[\|w_t^k - w_{t-1}^k\|_{U_k^{-1}}^2\big] \\
  &+ \mu_F^2u_k^{\min}(u_k^{\max} - \frac{2}{L_{\Omega}}) \mathbb{E}\big[\|w_t^k - w_{t-1}^k\|_{U_k^{-1}}^2\big] \\
   = & \big(1- \mu_F^2u_k^{\min}(\frac{2}{L_{\Omega}} - u_k^{\max})\big)\mathbb{E}\big[\|w_t^k - w_{t-1}^k\|_{U_k^{-1}}^2\big].
\end{align*}
Here, the first inequality holds due to the definition of $u_k^{\max}$, and the second inequality uses $\mathbb{E}\big[\|\nabla F(w_t^k) - \nabla F(w_{t-1}^k)\|_2^2\big] = \mathbb{E}\big[\|\mathbb{E}[v_t^k - v_{t-1}^k]\|_2^2\big] \leq \mathbb{E}\big[\|v_t^k - v_{t-1}^k\|_2^2\big]$, because it holds that $\mathbb{E}[\|z-\mathbb{E}[z]\|_2^2] = \mathbb{E}[\|z\|_2^2] - \|\mathbb{E}[z]\|^2 \geq 0$ for random vector $z\in \mathbb{R}^d$. Notice that $u_k^{\max} - 2/L_{\Omega}\leq0$ since $U_k\preceq 2/L_{\Omega}I$. In the third inequality we use the fact that $\mu_F\|w_t^k - w_{t-1}^k\|_2 \leq \|\nabla F(w_t^k) - \nabla F(w_{t-1}^k)\|_2$, which can be deduced from the strong convexity of $F(w)$. The last inequality is due to the definition of $u_k^{\min}$. The proof of the desired result is completed.
\end{IEEEproof}

The following theorem establishes the linear convergence of VM-mSRGBB under the strongly convex condition.
\begin{thm}\label{vmthm1}
Suppose that Assumptions \ref{vmnonsmooth}, \ref{vmsmooth}, and \ref{vmstrong} hold, and choose $b\in\{1,2,\ldots,n\}$. Assume that $0 \prec U_k\preceq 1/L_{\Omega}I$, $8L_{\Omega}u_k^{\max}/b<1$, and $m$ is chosen so that
\begin{equation*}
  \rho_k = \frac{1}{m\mu u_k^{\min}\big(1-\frac{8L_{\Omega}u_k^{\max}}{b}\big)} + \frac{4L_{\Omega}u_k^{\max}}{mb\big(1-\frac{8L_{\Omega}u_k^{\max}}{b}\big)} <1.
\end{equation*}
Then, VM-mSRGBB converges linearly in expectation:
\begin{equation*}
  \mathbb{E}\big[P(\tilde{w}^{k+1}) - P(w_*)\big] \leq \rho_k \mathbb{E}\big[P(\tilde{w}^{k}) - P(w_*)\big].
\end{equation*}
\end{thm}

\begin{IEEEproof}
From the update rule \eqref{eqvmspm2}, we obtain that, for any $t\geq1$,
\begin{align}\label{eqvm2.13}
    & \|w_{t+1}^k - w_*\|_{U_k^{-1}}^2\nonumber\\
    =& \|w_t^k - U_kg_t^k - w_*\|_{U_k^{-1}}^2 \nonumber\\
    =& \|w_t^k - w_*\|_{U_k^{-1}}^2 - 2(w_t^k- w_*)^Tg_t^k + \|g_t^k\|_{U_k}^2\nonumber\\
    \leq& \|w_t^k - w_*\|_{U_k^{-1}}^2 - 2\big(P(w_{t+1}^k) - P(w_*)\big) \nonumber\\
     & + 2(w_{t+1}^k - w_*)^T\delta_t^k,
\end{align}
where the last inequality uses Lemma \ref{vmlem5}. In order to provide an upper bound on the quantity $2(w_{t+1}^k - w_*)^T\delta_t^k$, we need the following notation
\begin{equation}\label{eqvmpgd}
  \bar{w}_{t+1}^k = \textrm{prox}_{R}^{U_k^{-1}}(w_t^k - U_k\nabla F(w_t^k)),
\end{equation}
which is independent of the random variable $I_t$. Then we get
\begin{align}\label{eqvmc4.1}
  &  2(w_{t+1}^k - w_*)^T\delta_t^k \nonumber\\
  = & 2(w_{t+1}^k - \bar{w}_{t+1}^k)^T\delta_t^k + 2(\bar{w}_{t+1}^k - w_*)^T\delta_t^k \nonumber\\
  \leq& 2\|\delta_t^k\|_{U_k}\|w_{t+1}^k - \bar{w}_{t+1}^k\|_{U_k^{-1}} + 2(\bar{w}_{t+1}^k - w_*)^T\delta_t^k \nonumber\\
  \leq& 2\|\delta_t^k\|_{U_k}\|(w_t^k - U_k v_t^{k})-(w_t^k - U_k \nabla F(w_t^k))\|_{U_k^{-1}} \nonumber\\
  & + 2(\bar{w}_{t+1}^k - w_*)^T\delta_t^k \nonumber\\
  \leq& 2u_k^{\max}\|\delta_t^k\|_2^2+ 2(\bar{w}_{t+1}^k - w_*)^T\delta_t^k,
\end{align}
where the first equality uses the fact that $|w^Tw'|^2 \leq \|w\|_A^2\cdot \|w'\|_{A^{-1}}^2$ with any symmetric positive definite matrix $A$, the second inequality holds due to Lemma \ref{vmlem6}, and the last inequality follows from the definition of $u_k^{\max}$ and $\delta_t^k$. Combining \eqref{eqvmc4.1} with \eqref{eqvm2.13}, we obtain
\begin{align}\label{eqvmc4.2}
  & \|w_{t+1}^k -w_*\|_{U_k^{-1}}^2 \nonumber\\
 \leq & \|w_t^k - w_*\|_{U_k^{-1}}^2 - 2\big(P(w_{t+1}^k) - P(w_*)\big)\nonumber\\
 & + 2u_k^{\max}\|\delta_t^k\|_2^2+ 2(\bar{w}_{t+1}^k - w_*)^T\delta_t^k.
\end{align}
Since both $\bar{w}_{t+1}^k$ and $w_*$ are independent of $I_t$ and the history of random variables $w_0^k$, $I_1$, $I_2$, $\ldots$, $I_{t-1}$, and $\mathbb{E}[\delta_t^k] = \mathbb{E}[\mathbb{E}[\nabla F(w_t^k)-v_t^{k} |\mathcal{F}_t]] = 0$, we have
\begin{equation*}\label{eqvm2.13.4}
  \mathbb{E}\big[(\bar{w}_{t+1}^k - w_*)^T\delta_t^k\big] = 0.
\end{equation*}
By taking expectation with respect to all the variables generated in the $k$-th outer loop and applying Lemma \ref{vmlem1} to \eqref{eqvmc4.2}, we obtain
\begin{align}\label{eqvm3.6}
  &  \mathbb{E}\big[\|w_{t+1}^k -w_*\|_{U_k^{-1}}^2\big]\nonumber\\
  \leq& \mathbb{E}\big[\|w_t^k - w_*\|_{U_k^{-1}}^2\big] - 2\mathbb{E}\big[P(w_{t+1}^k) - P(w_*)\big] \nonumber\\
  & + 2u_k^{\max}\mathbb{E}\big[\|\delta_t^k\|_2^2\big] \nonumber\\
  \leq& \mathbb{E}\big[\|w_t^k - w_*\|_{U_k^{-1}}^2\big] - 2\mathbb{E}\big[P(w_{t+1}^k) - P(w_*)\big] \nonumber\\
  & + \frac{8L_{\Omega}u_k^{\max}}{b}\mathbb{E}\big[P(w_t^k) - P(w_*)\big]\nonumber\\
  & + \frac{8L_{\Omega}u_k^{\max}}{b}\mathbb{E}\big[P(w_{t-1}^k) - P(w_*)\big].
\end{align}
Notice  that  $v_1^k = v_0^k$ and $\delta_1^k=\nabla F(w_1^k)-v_1^k=\nabla F(\tilde{w}^k)-v_0^k=0$ since $w_1^k=w_0^k=\tilde{w}^k$ and $v_0^k = \nabla F(\tilde{w}^k)$.  So, it follows from \eqref{eqvm2.13} that
\begin{equation}\label{eqvm3.62}
\|w_{2}^k - w_*\|_{U_k^{-1}}^2 \leq \|w_{1}^k - w_*\|_{U_k^{-1}}^2 - 2\big((P(w_{2}^k) - P(w_*)\big).
\end{equation}
Summing \eqref{eqvm3.6} over $t = 2, \ldots, m$ and taking into account \eqref{eqvm3.62}, we get
\begin{align}\label{eqvm3.7}
  &  \mathbb{E}\big[\|w_{m+1}^k - w_*\|_{U_k^{-1}}^2\big] + 2\mathbb{E} \big[P(w_{m+1}^k) - P(w_*)\big] \nonumber\\
  &+2\big(1-\frac{4L_{\Omega}u_k^{\max}}{b}\big)\sum_{t=2}^{m} \mathbb{E}\big[P(w_t^k) - P(w_*)\big]\nonumber\\
  \leq & \mathbb{E}\big[\|w_1^k - w_*\|_{U_k^{-1}}^2\big] + \frac{8L_{\Omega}u_k^{\max}}{b}\mathbb{E}\big[P(w_1^k) - P(w_*)\big] \nonumber\\
  & + \frac{8L_{\Omega}u_k^{\max}}{b}\sum_{t=2}^{m-1}\mathbb{E}\big[P(w_t^k) - P(w_*)\big] \nonumber\\
  \leq & \mathbb{E}\big[\|w_1^k - w_*\|_{U_k^{-1}}^2\big] + \frac{8L_{\Omega}u_k^{\max}}{b}\mathbb{E}\big[P(w_1^k) - P(w_*)\big] \nonumber\\
  & + \frac{8L_{\Omega}u_k^{\max}}{b}\sum_{t=2}^{m} \mathbb{E}\big[P(w_t^k) - P(w_*)\big],
\end{align}
where the last inequality uses the fact that $P(w_t^k) \geq P(w_*)$ for all $t\geq 0$. By rearranging terms of \eqref{eqvm3.7}, we get
\begin{align}\label{eqvm3.7.2}
  &  \mathbb{E}\big[\|w_{m+1}^k - w_*\|_{U_k^{-1}}^2\big] + 2\mathbb{E} \big[P(w_{m+1}^k) - P(w_*)\big] \nonumber\\
  &+2\big(1-\frac{8L_{\Omega}u_k^{\max}}{b}\big)\sum_{t=2}^{m} \mathbb{E}\big[P(w_t^k) - P(w_*)\big]\nonumber\\
  \leq & \mathbb{E}\big[\|w_1^k - w_*\|_{U_k^{-1}}^2\big] + \frac{8L_{\Omega}u_k^{\max}}{b}\mathbb{E}\big[P(w_1^k) - P(w_*)\big],
\end{align}
Since $2(1-\frac{8L_{\Omega}u_k^{\max}}{b}) < 2$, $\mathbb{E}\big[\|w_{m+1}^k - w_*\|_{U_k^{-1}}^2\big] \geq 0$, and $w_1^k = \tilde{w}^k$, we obtain
\begin{align*}\label{eqvm3.7.3}
  & 2\big(1-\frac{8L_{\Omega}u_k^{\max}}{b} \big) \sum_{t=2}^{m+1} \mathbb{E}\big[P(w_t^k) - P(w_*)\big] \nonumber\\
  \leq&  \|\tilde{w}^k - w_*\|_{U_k^{-1}}^2 + \frac{8L_{\Omega}u_k^{\max}}{b}\mathbb{E}\big[P(\tilde{w}^k) - P(w_*)\big] \nonumber\\
  \leq&  \frac{1}{u_k^{\min}}\cdot\|\tilde{w}^k - w_*\|_{2}^2 + \frac{8L_{\Omega}u_k^{\max}}{b}\mathbb{E}\big[P(\tilde{w}^k) - P(w_*)\big] \nonumber\\
  \leq& \left(\frac{2}{\mu u_k^{\min}} + \frac{8L_{\Omega}u_k^{\max}}{b}\right) \mathbb{E}\big[P(\tilde{w}^k) - P(w_*)\big],
\end{align*}
where the second inequality holds by the definition of $u_k^{\min}$ and in the last inequality we use the fact that $\|\tilde{w}^{k} - w_*\|_2^2 \leq \frac{2}{\mu}\big[P(\tilde{w}^{k}) - P(w_*)\big]$, which can be deduced from the strong convexity of $P(w)$. By the definition of $\tilde{w}^{k+1}$ in Algorithm \ref{conVM}, we have $\mathbb{E}[P(\tilde{w}^{k+1})] = \frac{1}{m}\sum_{t=1}^{m}\mathbb{E}[P(w_{t+1}^k)]$. Then the following inequality holds
\begin{align*}
  & 2m\big(1-\frac{8L_{\Omega}u_k^{\max}}{b}\big)\ \mathbb{E}\big[P(\tilde{w}^{k+1}) - P(w_*)\big] \\
   \leq& \left(\frac{2}{\mu u_k^{\min}} + \frac{8L_{\Omega}u_k^{\max}}{b}\right) \mathbb{E}\big[P(\tilde{w}^k) - P(w_*)\big].
\end{align*}
Dividing both sides of the above inequality by $2m\big(1-\frac{8L_{\Omega}u_k^{\max}}{b}\big)$ and using the definition of $\rho_k$, we arrive at
\begin{align*}
  & \mathbb{E}\big[P(\tilde{w}^{k+1}) - P(w_*)\big] \leq \rho_k \mathbb{E}\big[P(\tilde{w}^{k}) - P(w_*)\big].
\end{align*}
%\begin{align*}
%% \nonumber to remove numbering (before each equation)
%  & \mathbb{E}\big[P(\tilde{w}^{k+1}) - P(w_*)\big] \\
%  \leq & \Bigg(\frac{1}{m\mu u_k^{\min}\big(1-\frac{8L_{\Omega}u_k^{\max}}{b}\big)} + \frac{4L_{\Omega}u_k^{\max}}{mb\big(1-\frac{8L_{\Omega}u_k^{\max}}{b}\big)}\Bigg)\\
%  & \mathbb{E}\big[P(\tilde{w}^{k}) - P(w_*)\big]\\
%  = & \rho_k \mathbb{E}\big[P(\tilde{w}^{k}) - P(w_*)\big].
%\end{align*}
Then the desired result is proved.

\end{IEEEproof}

\subsection{VM-mSRGBB for non-strongly convex functions}
We establish linear convergence of our VM-mSRGBB method under quadratic growth condition (QGC) \cite{karimi2016linear}, which is stated as follows:
\begin{equation}\label{eqvmssc}
  P(w) - P_* \geq \frac{\nu}{2}\|w-\hat{w}\|_2^2, \textrm{ $\forall w\in \mathbb{R}^d$},
\end{equation}
where $\nu>0$, $\hat{w}$ is the projection of $w$ onto $\mathcal{W}_*$ and $P_*$ represents the optimal value of \eqref{eqvm1.1}.

QGC is weaker than the strongly convex condition. For example, the $\ell_1$-regularized least squares problems and logistic regression problems satisfying QGC \cite{gong2014linear}, however, they are not strongly convex when the data matrix does not have full column rank. It is shown that a nonsmooth convex function satisfies QGC meets the proximal Polyak-{\L}ojasiewicz inequality \cite{karimi2016linear}. The authors of \cite{zhang2017the} deduced the equivalence among QGC, the extended restricted strongly convex property (eRSC) and the extended global error bound property (eQEB).

\begin{thm}\label{vmthm2}
Suppose that Assumptions \ref{vmnonsmooth} and \ref{vmsmooth} hold, problem \eqref{eqvm1.1} satisfies QGC inequality with $\nu >0$, and choose $b\in\{1,2,\ldots,n\}$. Further assume that $0 \prec U_k\preceq 1/L_{\Omega}I$, $8L_{\Omega}u_k^{\max}/b<1$, and $m$ is chosen so that
\begin{equation*}
  \hat{\rho}_k = \frac{1}{m\nu u_k^{\min}\big(1-\frac{8L_{\Omega}u_k^{\max}}{b}\big)} + \frac{4L_{\Omega}u_k^{\max}}{mb\big(1-\frac{8L_{\Omega}u_k^{\max}}{b}\big)} <1.
\end{equation*}
Then, VM-mSRGBB achieves a linear convergence rate in expectation:
\begin{equation*}
  \mathbb{E}\big[P(\tilde{w}^{k+1}) - P_*\big] \leq \hat{\rho}_k \mathbb{E}\big[P(\tilde{w}^{k}) - P_*\big].
\end{equation*}
\end{thm}

\begin{IEEEproof}
Let $\hat{w}_{t}^k$ be the projection of $w_{t}^k$ onto $\mathcal{W}_*$, i.e.,
 $$\hat{w}_{t}^k = \Pi_{\mathcal{W}_*}(w_{t}^k) = \arg\min_w\{w\in \mathcal{W}_*: \|w_{t}^k - w\|_{U_k^{-1}}^2\}.$$
 Then $\hat{w}_{t}^k,\hat{w}_{t+1}^k\in \mathcal{W}_*$, which together with  \eqref{eqvmspm2} implies that, for $t\geq1$,
\begin{align*}
  &\|w_{t+1}^k - \hat{w}_{t+1}^k\|_{U_k^{-1}}^2 \nonumber\\
  \leq& \|w_{t+1}^k - \hat{w}_{t}^k\|_{U_k^{-1}}^2 \nonumber\\
   =& \|w_t^k -U_kg_t^k- \hat{w}_{t}^k\|_{U_k^{-1}}^2 \nonumber\\
   =&\|w_t^k -\hat{w}_{t}^k\|_{U_k^{-1}}^2 + 2\eta_k(\hat{w}_{t}^k -w_t^k)^Tg_t^k + \|g_t^k\|_{U_k}^2 \nonumber\\
   \leq & \|w_t^k -\hat{w}_{t}^k\|_{U_k^{-1}}^2 + 2(P_* - P(w_{t+1}^k)) - 2(\hat{w}_t^k-w_{t+1}^k)^T\delta_t^k,\nonumber
\end{align*}
where the first inequality holds due to the positive definiteness of $U_k$, and the last inequality is the application of Lemma \ref{vmlem5} with $\hat{w}_{t}^k\in \mathcal{W}_*$.

%Similar to the proving process of \eqref{eqvmc4.1}, by using the proximal full gradient in \eqref{eqvmpgd}, we have
%\begin{equation}\label{eqvm5.1}
%  - 2(\hat{w}_t^k-w_{t+1}^k)^T\delta_t^k \leq 2u_k^{\max}\|\delta_t^k\|_2^2+ 2(\bar{w}_{t+1}^k - \hat{w}_t^k)^T\delta_t^k.
%\end{equation}
Similarly to the proof of \eqref{eqvm3.62}-\eqref{eqvm3.7.2} in Theorem \ref{vmthm1}, we obtain
\begin{align}\label{eqvmc8}
  & 2\big(1-\frac{8L_{\Omega}u_k^{\max}}{b} \big) \sum_{t=2}^{m+1} \mathbb{E}\big[P(w_t^k) - P_*\big] \nonumber\\
  \leq&  \|\tilde{w}^k - \hat{w}_1^k\|_{U_k^{-1}}^2 + \frac{8L_{\Omega}u_k^{\max}}{b}\mathbb{E}\big[P(\tilde{w}^k) - P_*\big]\nonumber\\
  \leq & \frac{1}{u_k^{\min}}\cdot\|\tilde{w}^k - \hat{w}_1^k\|_2^2 + \frac{8L_{\Omega}u_k^{\max}}{b}\mathbb{E}\big[P(\tilde{w}^k) - P_*\big]
\end{align}
The definition of $\tilde{w}^{k+1}$ implies that $\mathbb{E}[P(\tilde{w}^{k+1})] =\frac{1}{m}\sum_{t=1}^{m} \mathbb{E}[P(w_{t+1}^k)]$. Considering QGC with $w = \tilde{w}^k$, $\tilde{w}^k =w_1^k$ and $\hat{w}_{1}^k = \Pi_{\mathcal{W_*}} (w_1^k) \in \mathcal{W}_*$, we have
\begin{equation*}
  P(\tilde{w}^k) - P_* \geq \frac{\nu}{2}\| \tilde{w}^k - \hat{w}_{1}^k\|_2^2,
\end{equation*}
which together with \eqref{eqvmc8} yields
\begin{align*}
  &2m\big(1-\frac{8L_{\Omega}u_k^{\max}}{b}\big)\mathbb{E}\big[P(\tilde{w}^{k+1}) - P_*\big]\\
  \leq & \left(\frac{2}{\nu u_k^{\min}} + \frac{8L_{\Omega}u_k^{\max}}{b}\right) \mathbb{E}\big[P(\tilde{w}^k) - P_*\big].
\end{align*}
Dividing both sides of the above inequality by $2m(1-\frac{8L_{\Omega}u_k^{\max}}{b})$, and considering the definition of $\hat{\rho}_k$, we arrive at
\begin{equation*}
   \mathbb{E}\big[P(\tilde{w}^{k+1}) - P_*\big]
  \leq  \hat{\rho}_k\mathbb{E}\big[P(\tilde{w}^{k}) - P_*\big].
\end{equation*}
\end{IEEEproof}

\subsection{VM-mSRGBB for convex functions}
Now we study the convergence of VM-mSRGBB for convex nonsmooth functions.

Next lemma presents a new 3-point property which generalizes the one in \cite{lan2012optimal}.
\begin{lem}(generalized 3-point property)\label{vmlem7}
Suppose that $R:\mathbb{R}^d\rightarrow \mathbb{R}$ is lower semicontinuous convex (but possibly nondifferentiable) and $w' = \textrm{prox}_R^{A^{-1}}(w)$ with $A\in\mathbb{S}_{++}^{d \times d}$. Then, for any $z\in \mathbb{R}^d$, we have the following inequality:
\begin{equation*}
  R(w') + \frac{1}{2}\|w'-w\|_{A^{-1}}^2 \leq R(z) + \frac{1}{2}\|z-w\|_{A^{-1}}^2 - \frac{1}{2}\|w'-z\|_{A^{-1}}^2.
\end{equation*}
\end{lem}
\begin{IEEEproof}
Since $w' = \textrm{prox}_R^{A^{-1}}(w) = \arg\min_z\{R(z) + \frac{1}{2}\|z-w\|_{A^{-1}}^2\}$, there exists $\varpi \in \partial R(w')$ such that
\begin{equation*}
  \varpi + A^{-1}(w'-w)= 0.
\end{equation*}
By direct expansion, we have
\begin{align*}
  \frac{1}{2}\|z-w\|_{A^{-1}}^2 =& \frac{1}{2}\|z-w'\|_{A^{-1}}^2 + \frac{1}{2}\|w'-w\|_{A^{-1}}^2 \\
   & + (z-w')^TA^{-1}(w'-w), \quad \forall z\in\mathbb{R}^d.
\end{align*}
Using the above two relations and the convexity of  $R(z)$, we conclude that
\begin{align*}
  & R(z) + \frac{1}{2}\|z-w\|_{A^{-1}}^2\\
   =&  R(z) +\frac{1}{2}\|z-w'\|_{A^{-1}}^2 + \frac{1}{2}\|w'-w\|_{A^{-1}}^2 \\
   &+ (z-w')^TA^{-1}(w'-w) \\
   \geq& R(w') + \varpi^T(z-w')  + \frac{1}{2}\|z-w'\|_{A^{-1}}^2 + \frac{1}{2}\|w'-w\|_{A^{-1}}^2 \\
   & + (z-w')^TA^{-1}(w'-w)\\
   = & R(w') + \frac{1}{2}\|z-w'\|_{A^{-1}}^2+ \frac{1}{2}\|w'-w\|_{A^{-1}}^2.
\end{align*}
\end{IEEEproof}

\begin{lem}\label{vmlem8}
Suppose that $R:\mathbb{R}^d\rightarrow \mathbb{R}$ is lower semicontinuous convex (but possibly nondifferentiable) and
\begin{equation}\label{eqvm6.2}
  w' = \textrm{prox}_R^{A^{-1}}(w-A \zeta)
\end{equation}
with $A\in\mathbb{S}_{++}^{d \times d}$ and $\zeta\in \mathbb{R}^d$. Then, the following inequality holds
\begin{align}\label{eqvm6.3}
  R(w')\leq& R(z)+ (z-w')^T\zeta\nonumber\\
& + \frac{1}{2}\big[\|z-w\|_{A^{-1}}^2 - \|w'-w\|_{A^{-1}}^2 - \|w'-z\|_{A^{-1}}^2\big]
\end{align}
for all $z\in \mathbb{R}^d$.
\end{lem}
\begin{IEEEproof}
By applying Lemma \ref{vmlem7} to \eqref{eqvm6.2}, we get
\begin{align}\label{eqvm6.4}
  & R(w') + (w'-w)^T\zeta+ \frac{1}{2}\|w'-w\|_{A^{-1}}^2 + \frac{1}{2} \|\zeta\|_{A}^2\nonumber\\
  =& R(w') + \frac{1}{2}\|w'-(w -A\zeta)\|_{A^{-1}}^2 \nonumber\\
  \leq & R(z) + \frac{1}{2}\|z-(w-A\zeta)\|_{A^{-1}}^2 - \frac{1}{2}\|w'-z\|_{A^{-1}}^2\nonumber\\
  =& R(z) + (z-w)^T\zeta + \frac{1}{2}\|z-w\|_{A^{-1}}^2 + \frac{1}{2} \|\zeta\|_{A}^2 \nonumber\\
  & - \frac{1}{2}\|w'-z\|_{A^{-1}}^2.
\end{align}
\end{IEEEproof}

\begin{lem}\label{vmlem9}
Consider $P(w)$ as defined in \eqref{eqvm1.1}. Suppose that Assumptions \ref{vmnonsmooth} and \ref{vmsmooth} holds. Then, for $w'$ defined by \eqref{eqvm6.2}, the following inequality holds:
\begin{align*}
  P(w') \leq& P(z) + (w'-z)^T(\nabla F(w) - \zeta) \\
   & + \frac{1}{2}\|w'-w\|_{(L_{\Omega}I-A^{-1})}^2 + \frac{1}{2}\|z-w\|_{(L_{\Omega}I+A^{-1})}^2 \\
   &- \frac{1}{2}\|w'-z\|_{A^{-1}}^2,
\end{align*}
for all $z\in \mathbb{R}^d$.
\end{lem}
\begin{IEEEproof}
From the Lipschitz continuity of $\nabla F$ and the fact $L\geq L_{\Omega}$, we obtain
\begin{align*}
  F(w') \leq& F(w) + \nabla F(w)^T(w'-w) + \frac{L_{\Omega}}{2}\|w'-w\|_2^2, \\
   F(w) \leq&  F(z) + \nabla F(w)^T(w-z) + \frac{L_{\Omega}}{2}\|w-z\|_2^2.
\end{align*}
By summing the above two inequalities, we have
\begin{align}\label{eqvm6.5}
  F(w') \leq & F(z) + \nabla F(w)^T(w'-z) \nonumber\\
  & + \frac{L_{\Omega}}{2}\|w'-w\|_2^2+ \frac{L_{\Omega}}{2}\|w-z\|_2^2.
\end{align}
Summing \eqref{eqvm6.3} and \eqref{eqvm6.5}, we get
\begin{align*}
  P(w') \leq& P(z) + (w'-z)^T(\nabla F(w)-\zeta) \\
   & + \frac{1}{2}\|w'-w\|_{(L_{\Omega}I-A^{-1})}^2 + \frac{1}{2}\|z-w\|_{(L_{\Omega}I+A^{-1})}^2 \\
   & - \frac{1}{2}\|w'-z\|_{A^{-1}}^2,
\end{align*}
which completes our proof.
\end{IEEEproof}

In order to derive an upper bound on the variance of $v_t^k$ in the mini-batch setting, we first show the result in the case where $b=1$.
\begin{lem}\label{vmlem10.1}
Suppose that Assumption \ref{vmnonsmooth} holds. Consider $v_t^k$ as defined in \eqref{eqvmvt1} with $b=1$, i.e.,
\begin{equation}\label{eqvmvt1.1}
  v_t^k = \frac{\nabla f_{i_t}(w_t^k) - \nabla f_{i_t}(w_{t-1}^k)}{nq_{i_t}} + v_{t-1}^k,
\end{equation}
then the following inequality holds:
\begin{equation*}
  \mathbb{E}[\|v_t^k-\nabla F(w_t^k)\|_2^2] \leq L_{\Omega}^2\mathbb{E}[\|w_t^k - w_{t-1}^k\|_2^2], \quad \forall t\geq1.
\end{equation*}
\end{lem}
\begin{IEEEproof}
See Appendix \ref{vmapp2}.
\end{IEEEproof}

The following lemma provides an upper bound on $v_t^k$, which looks similar to the Lemma 3 of \cite{reddi2016proximal}, but they are essentially different due to the update rule of $v_t^k$.
\begin{lem}\label{vmlem10}
Suppose that Assumption \ref{vmnonsmooth} holds and choose $b\in\{1,2,\ldots,n\}$.  Consider $v_t^k$ as defined in \eqref{eqvmvt1}. Then, for any $t\geq1$, the following inequality holds
\begin{equation*}
  \mathbb{E}[\|v_t^k-\nabla F(w_t^k)\|_2^2] \leq \frac{L_{\Omega}^2}{b}\mathbb{E}[\|w_t^k - w_{t-1}^k\|_2^2].
\end{equation*}
\end{lem}
\begin{IEEEproof}
See Appendix \ref{vmapp3}.
\end{IEEEproof}

To establish the convergence of VM-mSRGBB under convex condition, we need the following notation of gradient mapping
\begin{equation}\label{eqvm6.18}
\mathcal{G}_{A^{-1}}(w) = A^{-1}\Big(w-\textrm{prox}_R^{A^{-1}} \big(w- A\nabla F(w)\big)\Big),
\end{equation}
where $A$ is a symmetric positive definite matrix. Note that when $R(w)$ is a constant function, the gradient mapping reduces to $\mathcal{G}_{A^{-1}}(w) = \nabla F(w)$. It is not difficult to show that $\mathcal{G}_{A^{-1}}(w)=0$ if and only if $w$ is a solution of problem \eqref{eqvm1.1}.

\begin{thm}\label{thm3}
Suppose that Assumptions \ref{vmnonsmooth} and \ref{vmsmooth} hold, and $0 \prec U_k\preceq 1/(3L_{\Omega})I$. Let $c_{t_k+1}=0$ and $c_t^k = c_{t+1}^k + \frac{u_k^{\max}L_{\Omega}^2}{2b}$. Then, for the output $w_a$ of Algorithm \ref{conVM}, after $T$ iterations, we have
\begin{equation*}
  \mathbb{E}[\|\mathcal{G}_{U_k^{-1}}(w_a)\|_{U_k}^2] \leq \frac{6(P(\tilde{w}^0) - P(w_*))}{T},
\end{equation*}
where $T = \sum_{k=0}^{K-1}t_k$.
\end{thm}
\begin{IEEEproof}
By applying Lemma \ref{vmlem9} to the proximal full gradient update defined in \eqref{eqvmpgd} (with $w' = \bar{w}_{t+1}^k$, $w=z=w_t^k$, $A = U_k$ and $\zeta = \nabla F(w_t^k)$), and taking total expectation over the entire history in the $k$-th outer loop, we have
\begin{equation}\label{eqvm6.19}
  \mathbb{E}[P(\bar{w}_{t+1}^k)] \leq \mathbb{E}[P(w_t^k) + \|\bar{w}_{t+1}^k - w_t^k\|_{(\frac{L_{\Omega}}{2}I - U_k^{-1})}^2].
\end{equation}
Recalling that the iterates of Algorithm \ref{conVM} are computed by
\begin{equation*}\label{eqvm6.20}
  w_{t+1}^k = \textrm{prox}_{R}^{U_k^{-1}}(w_t^k- U_kv_t^{k}).
\end{equation*}
Again by applying Lemma \ref{vmlem9} to the above update equation (with $w'=w_{t+1}^k$, $z=\bar{w}_{t+1}^k$, $w=w_t^k$, $A=U_k$ and $\zeta = v_t^k$) and taking expectation, we have
\begin{align}\label{eqvm6.21}
  & \mathbb{E}[P(w_{t+1}^k)] \nonumber\\
  \leq&  \mathbb{E}[P(\bar{w}_{t+1}^k)+ \frac{1}{2}\|\bar{w}_{t+1}^k - w_t^k\|_{(L_{\Omega}I+U_k^{-1})}^2 \nonumber\\
  & + \frac{1}{2}\|w_{t+1}^k - w_t^k\|_{(L_{\Omega}I-U_k^{-1})}^2 - \frac{1}{2}\|w_{t+1}^k - \bar{w}_{t+1}^k\|_{U_k^{-1}}\nonumber\\
   &+ (w_{t+1}^k - \bar{w}_{t+1}^k)^T(\nabla F(w_t^k) - v_t^k)].
\end{align}
By summing \eqref{eqvm6.19} and \eqref{eqvm6.21}, we obtain
\begin{align}\label{eqvm6.22}
  &\mathbb{E}[P(w_{t+1}^k)] \nonumber\\
  \leq& \mathbb{E}[P(w_t^k) +  \|\bar{w}_{t+1}^k - w_t^k\|_{(L_{\Omega}I-\frac{1}{2}U_k^{-1})}^2 \nonumber\\
  &+ \frac{1}{2}\|w_{t+1}^k - w_t^k\|_{(L_{\Omega}I-U_k^{-1})}^2 - \frac{1}{2}\|w_{t+1}^k - \bar{w}_{t+1}^k\|_{U_k^{-1}} \nonumber\\
   & + (w_{t+1}^k - \bar{w}_{t+1}^k)^T(\nabla F(w_t^k) - v_t^k)].
\end{align}
Let $\Gamma= (w_{t+1}^k - \bar{w}_{t+1}^k)^T(\nabla F(w_t^k) - v_t^k)$. The expectation on $\Gamma$ can be bounded above by
\begin{align*}
  \mathbb{E}[\Gamma] &\leq \frac{1}{2}\mathbb{E}[\|w_{t+1}^k - \bar{w}_{t+1}^k\|_{U_k^{-1}}] + \frac{1}{2} \mathbb{E}[\|\nabla F(w_t^k) - v_t^k\|_{U_k}^2] \\
   &\leq \frac{1}{2}\mathbb{E}[\|w_{t+1}^k - \bar{w}_{t+1}^k\|_{U_k^{-1}}] + \frac{u_k^{\max}L_{\Omega}^2}{2b} \mathbb{E}[\|w_t^k - w_{t-1}^k\|_2^2],
\end{align*}
where in the first inequality we use Cauchy-Schwarz and Young's inequality, and the second inequality follows from the definition of $u_k^{\max}$ and Lemma \ref{vmlem10}. We substitute the upper bound on $\Gamma$ in \eqref{eqvm6.22} and then obtain
\begin{align}\label{eqvm6.23}
  &\mathbb{E}[P(w_{t+1}^k)] \nonumber\\
  \leq& \mathbb{E}[P(w_t^k) +  \|\bar{w}_{t+1}^k - w_t^k\|_{(L_{\Omega}I-\frac{1}{2}U_k^{-1})}^2 \nonumber\\
   &  + \frac{1}{2}\|w_{t+1}^k - w_t^k\|_{(L_{\Omega}I-U_k^{-1})}^2 + \frac{u_k^{\max}L_{\Omega}^2}{2b}\|w_t^k - w_{t-1}^k\|_2^2].
\end{align}
In order to further analyze \eqref{eqvm6.23}, we need the following auxiliary function:
\begin{equation}\label{eqvm6.24}
 \Upsilon(w_{t+1}^k) = \mathbb{E}[P(w_{t+1}^k) + c_{t+1}^k\|w_{t+1}^k - w_t^k\|_2^2],
\end{equation}
where $c_{t_k+1}^k = 0$, and $c_t^k = c_{t+1}^k + \frac{u_k^{\max}L_{\Omega}^2}{2b}$. Then $\Upsilon(w_{t+1}^k)$ can be bounded above by
\begin{align}\label{eqvm6.25}
  & \Upsilon(w_{t+1}^k) \nonumber\\
  = &  \mathbb{E}[P(w_{t+1}^k) + c_{t+1}^k\|w_{t+1}^k - w_t^k\|_2^2] \nonumber\\
  \leq & \mathbb{E}[P(w_{t+1}^k) + c_{t+1}^k\|w_t^k - w_{t-1}^k\|_2^2] \nonumber\\
  \leq & \mathbb{E}[P(w_t^k) +  \|\bar{w}_{t+1}^k - w_t^k\|_{(L_{\Omega}I-\frac{1}{2}U_k^{-1})}^2  \nonumber\\
  &+ (c_{t+1}^k+\frac{u_k^{\max}L_{\Omega}^2}{2b})\|w_t^k - w_{t-1}^k\|_2^2]\nonumber\\
  = & \Upsilon(w_t^k) + \mathbb{E}[\|\bar{w}_{t+1}^k - w_t^k\|_{(L_{\Omega}I-\frac{1}{2}U_k^{-1})}^2],
\end{align}
where the first inequality follows from Theorem \ref{vmthmdes}, the second inequality holds by \eqref{eqvm6.23} and $0 \prec U_k\preceq 1/(3L_{\Omega})I \prec 1/L_{\Omega}I$, and the last equality is due to the definitions of $c_t^k$ and $\Upsilon(w_t^k)$. By summing \eqref{eqvm6.25} over $t=1,\ldots,t_k$, we get
\begin{equation}\label{eqvm6.26}
  \Upsilon(w_{t_k+1}^k) \leq \Upsilon(w_1^k) + \sum_{t=1}^{t_k} \mathbb{E}[\|\bar{w}_{t+1}^k - w_t^k\|_{(L_{\Omega}I-\frac{1}{2}U_k^{-1})}^2].
\end{equation}
By the fact $c_{t_k+1}^k=0$ and the definition of $\tilde{w}^{k+1}$, we have
$$\Upsilon(w_{t_k+1}^k) = \mathbb{E}[P(w_{t_k+1}^k)] = \mathbb{E}[P(\tilde{w}^{k+1})].$$
 Since $w_1^k = w_0^k =\tilde{w}^k$, we know that $\Upsilon(w_1^k) = \mathbb{E}[P(w_1^k)] = \mathbb{E}[P(\tilde{w}^k)]$. It follows from \eqref{eqvm6.26} that
\begin{equation}\label{eqvm6.27}
  \mathbb{E}[P(\tilde{w}^{k+1})] \leq \mathbb{E}[P(\tilde{w}^k)]+ \sum_{t=1}^{t_k}  \mathbb{E}[\|\bar{w}_{t+1}^k - w_t^k\|_{(L_{\Omega}I-\frac{1}{2}U_k^{-1})}^2].
\end{equation}
By summing \eqref{eqvm6.27} over $k=0,\ldots,K-1$ and rearranging  terms, we obtain
\begin{align}\label{eqvm6.28}
  \sum_{k=0}^{K-1}\sum_{t=1}^{t_k} \mathbb{E}[\|\bar{w}_{t+1}^k - w_t^k\|_{(\frac{1}{2}U_k^{-1}-L_{\Omega}I)}^2] \leq& P(\tilde{w}^0) - P(\tilde{w}^K) \nonumber\\
  \leq& P(\tilde{w}^0) - P(w_*),
\end{align}
where in the second inequality we use the fact that $P(\tilde{w}^k) \geq P(w_*)$ for all $k\in \{0,1,\ldots,K\}$.

From \eqref{eqvm6.18} and \eqref{eqvmpgd}, it follows that
\begin{align*}
  \mathcal{G}_{U_k^{-1}}(w_t^k) =& U_k^{-1}\Big(w_t^k - \textrm{prox}_R^{U_k^{-1}} \big(w_t^k - U_k\nabla F(w_t^k)\big)\Big) \\
  =& U_k^{-1}\Big(w_t^k - \bar{w}_{t+1}^k\Big).
\end{align*}
By using the fact $0 \prec U_k \preceq 1/(3L_{\Omega})I$, we have
\begin{align*}\label{eqvm6.292}
&\|\bar{w}_{t+1}^k - w_t^k\|_{(\frac{1}{2}U_k^{-1}-L_{\Omega}I)}^2\\
=& \|U_k\mathcal{G}_{U_k^{-1}}\|_{(\frac{1}{2}U_k^{-1}-L_{\Omega}I)}^2\\
=&  \mathcal{G}_{U_k^{-1}}^TU_k^T (\frac{1}{2}U_k^{-1}-L_{\Omega}I)U_k\mathcal{G}_{U_k^{-1}} \\
\geq & \mathcal{G}_{U_k^{-1}}^TU_k^T(\frac{1}{6}U_k^{-1})U_k\mathcal{G}_{U_k^{-1}}\\
=& \frac{1}{6}\|\mathcal{G}_{U_k^{-1}}\|_{U_k}^2.
\end{align*}
Combining the above inequality with \eqref{eqvm6.28}, we get
\begin{equation}\label{eqvm6.29}
   \sum_{k=0}^{K-1}\sum_{t=1}^{t_k} \frac{1}{6}\mathbb{E}[\|\mathcal{G}_{U_k^{-1}}\|_{U_k}^2] \leq P(\tilde{w}^0) - P(w_*).
\end{equation}
Then we obtain the desired result by the definitions of $w_a$ and $T$.
\end{IEEEproof}

\section{Numerical experiments}\label{experiment}
In this section, we present experimental results on the following  elastic net regularized logistic regression problem
\begin{equation}\label{eqvm7.2}
\min_{w\in\mathbb{R}^d} P(w)=\frac{1}{n}\sum_{i=1}^n \log(1+\exp(-b_ia_i^Tw))+\frac{\lambda_2}{2}\|w\|_2^2 +\lambda_1\|w\|_1,
\end{equation}
which is usually employed in machine learning for binary classification. All the test were performed with $R(w) = \lambda_1\|w\|_1$ and
\begin{equation*}
f_i(w) = \log(1+\exp(-b_ia_i^Tw))+\frac{\lambda_2}{2}\|w\|_2^2.
\end{equation*}

Four publicly available data sets ijcnn1, rcv1, real-sim and covtype, which can be downloaded from the LIBSVM website \footnote{\url{www.csie.ntu.edu.tw/~cjlin/libsvmtools/}}, were tested. Table \ref{tab4} lists the detailed information of these four data sets, including their sizes $n$, dimensions $d$, and Lipschitz constants $L$. Moreover, the values of regularization parameters $\lambda_1$ and $\lambda_2$ used in our experiments are also listed in Table \ref{tab4}. Notice that the choices of regularization parameters are typical in machine learning benchmarks to obtain good classification performance, see \cite{xiao2014proximal} for example.
\begin{table}[!t]
\renewcommand{\arraystretch}{1.3}
\caption{Data sets and parameters used in numerical experiments}
\label{tab4}
\centering
\begin{tabular}{ccccccc}
  \hline
  Data sets & $n$ & $d$  & $\lambda_2$ & $\lambda_1$ & $L$\\
  \hline
  ijcnn1 & 49,990 & 22 &$10^{-4}$ & $10^{-5}$& 0.9842\\

  rcv1 & 20,242 & 47,236 &$10^{-4}$ & $10^{-5}$ & 0.2501\\

  real-sim & 72,309 & 20,958 &$10^{-4}$ & $10^{-5}$& 0.2501\\

  covtype & 581,012 & 54 &$10^{-5}$ & $10^{-4}$ & 1.9040\\

  \hline
\end{tabular}
\end{table}

For fair comparison, all methods were implemented in Matlab 2018b, and the experiments were conducted on a laptop with an Intel Core i7, 1.80 GHz processor and 16 GB of RAM running Windows 10 system. In Figs. 1-3, the $x$-axis is the number of effective passes over the data, where the evaluation of $n$ component gradients counts as one effective pass. The $y$-axis with ``optimality gap'' denotes the value $P(\tilde{w}^k) - P(w_*)$ with $w_*$ obtained by running Prox-SVRG with best-tuned fixed stepsizes.

\subsection{Comparison with Prox-SVRG and Prox-SVRG-BB}\label{ne2}

This subsection presents the results of VM-mSRGBB with $b=1$ for solving \eqref{eqvm7.2} on the four data sets listed in Table \ref{tab4}. Prox-SVRG and the proximal version of SVRG-BB  (Prox-SVRG-BB) were also run for comparison. Notice that the SVRG-BB method is proposed to solve problem \eqref{eqvm1.1} with $R(w)=0$. In order to solve the nonsmooth problem \eqref{eqvm7.2}, the proximal operator was incorporated to obtain the Prox-SVRG-BB method. For Prox-SVRG, as suggested in \cite{xiao2014proximal}, we set $m = 2n$. The best-tuned $m$ was employed by Prox-SVRG-BB.

\begin{figure}
%\centering
%\begin{minipage}[!t]{0.5\linewidth}
  \subfigure[ijcnn1]{
    \label{fig:subfig:ijcnn1a} %% label for first subfigure
    \includegraphics[width=1.6in,height=1.3in]{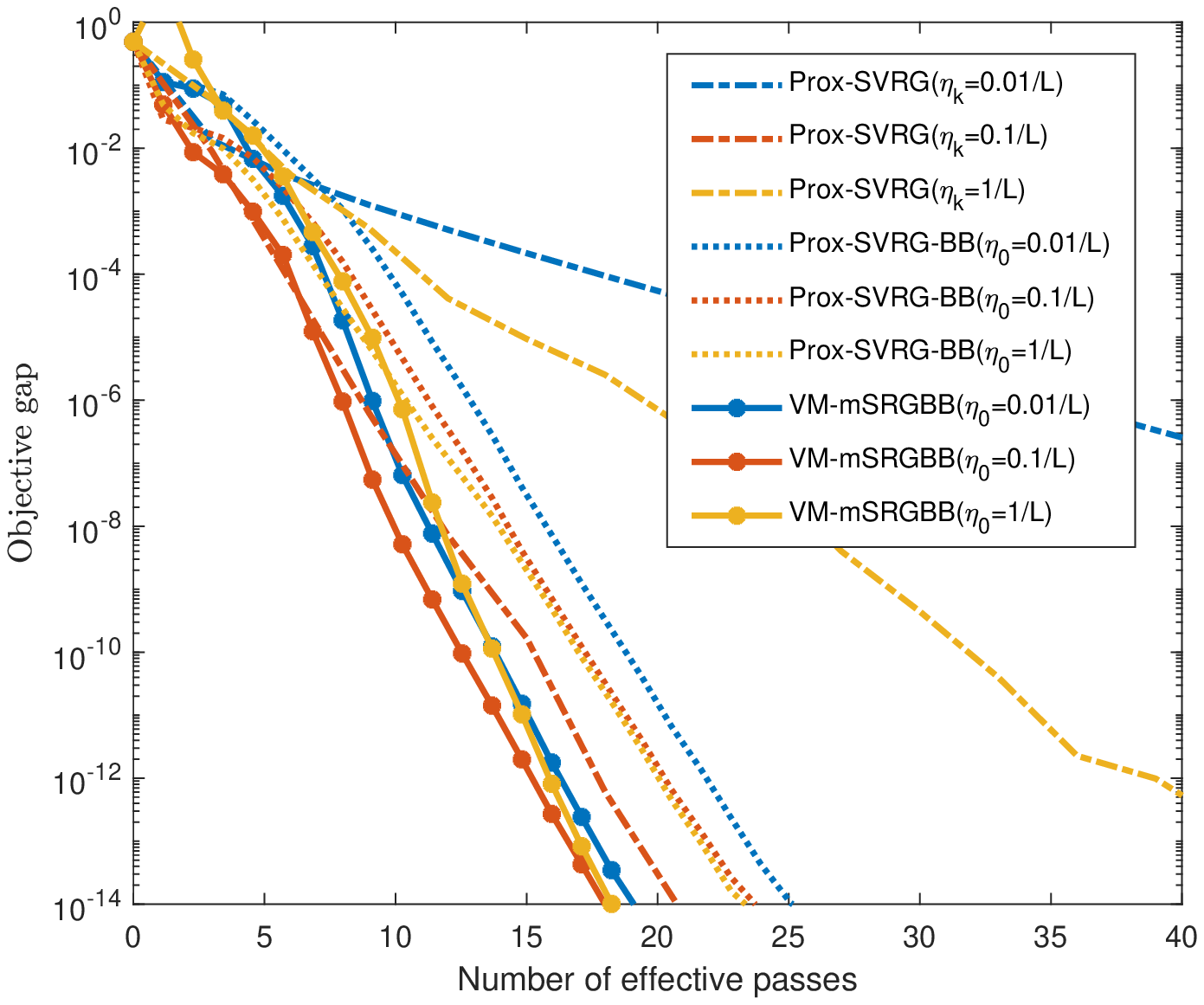}}
  \hspace{3mm}
  \subfigure[rcv1]{
    \label{fig:subfig:rcv1b} %% label for first subfigure
    \includegraphics[width=1.6in,height=1.3in]{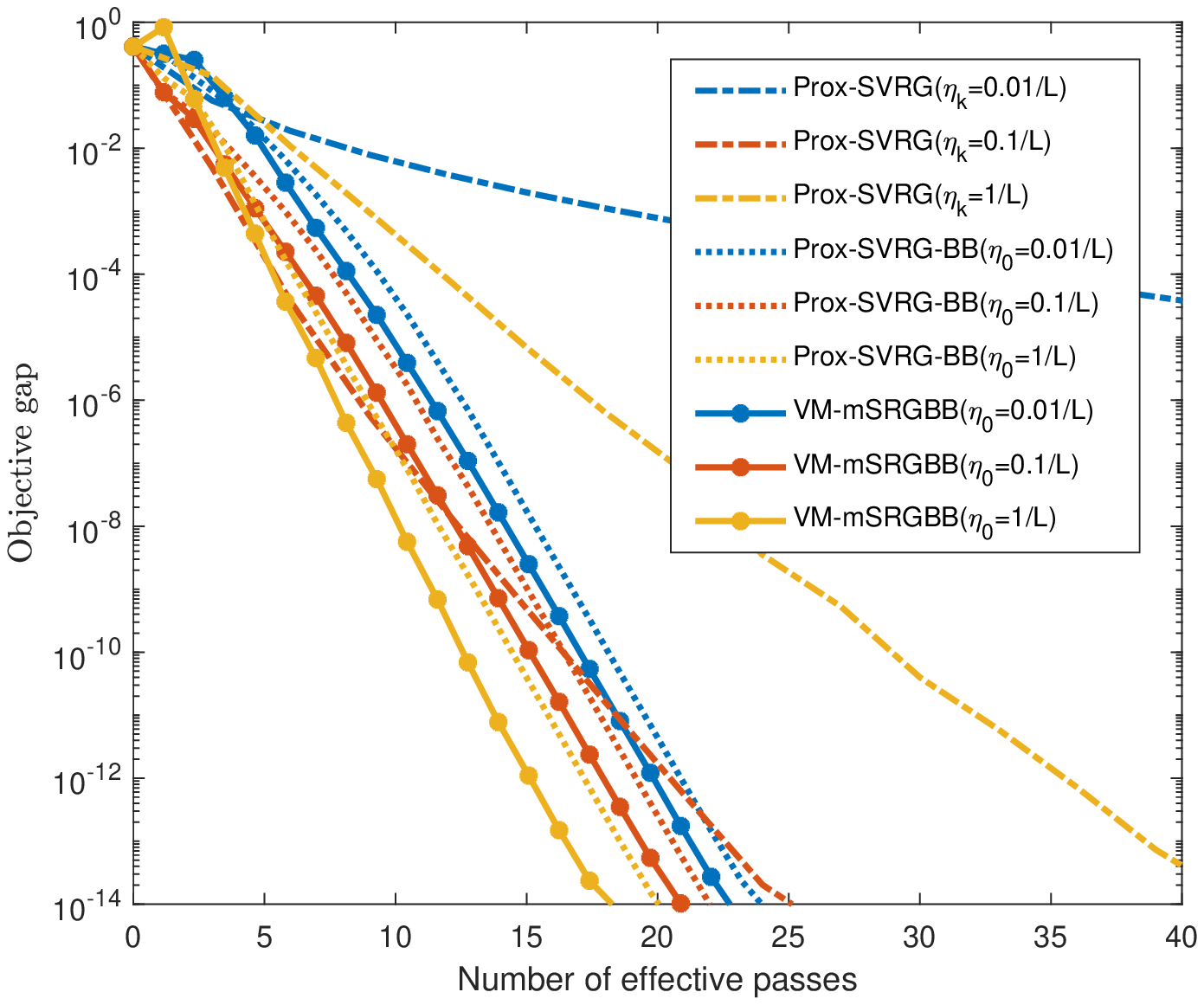}}
  \subfigure[real-sim]{
    \label{fig:subfig:real-simc} %% label for first subfigure
    \includegraphics[width=1.6in,height=1.3in]{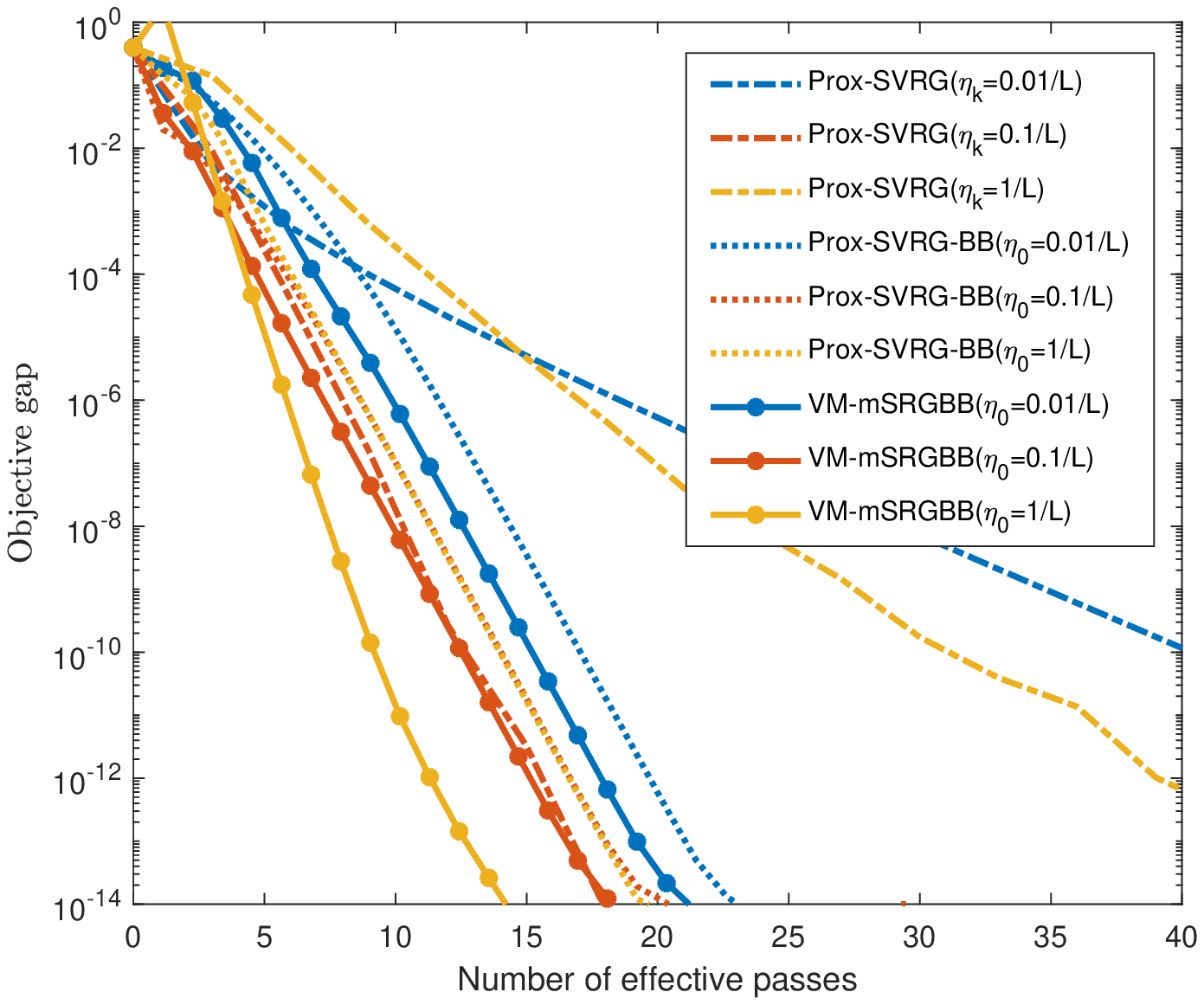}}
  \hspace{3mm}
  \subfigure[covtype]{
    \label{fig:subfig:covtyped} %% label for first subfigure
    \includegraphics[width=1.6in,height=1.3in]{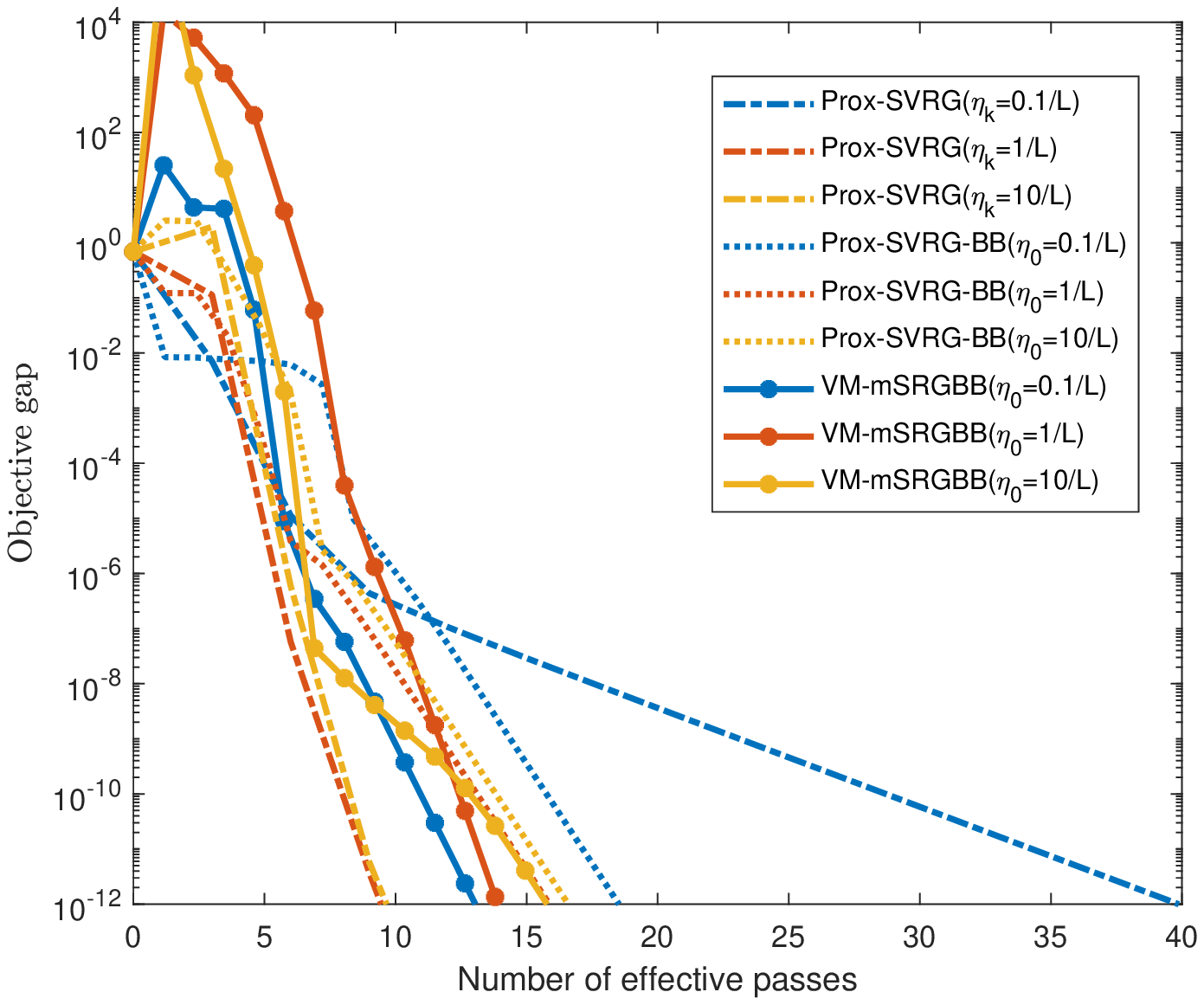}}
  \caption{Comparison of VM-mSRGBB, Prox-SVRG and Prox-SVRG-BB with different initial stepsizes.}\label{VSmS2GD}
\end{figure}

It can be seen from Fig. \ref{VSmS2GD} that VM-mSRGBB often performs better than Prox-SVRG with different initial stepsizes. Unlike Prox-SVRG, VM-mRGBB is not sensitive to the choice of initial stepsize, which would save much time on choosing initial stepsize so that it has promising potential in practice. Moreover, for different initial stepsizes, VM-mSRGBB performs better than Prox-SVRG-BB.

\subsection{Properties of VM-mSRGBB with different $b$}\label{ne1}
\begin{figure}
%\centering
%\begin{minipage}[!t]{0.5\linewidth}
  \subfigure[ijcnn1]{
    \label{fig:subfig:ijcnn1a} %% label for first subfigure
    \includegraphics[width=1.6in,height=1.3in]{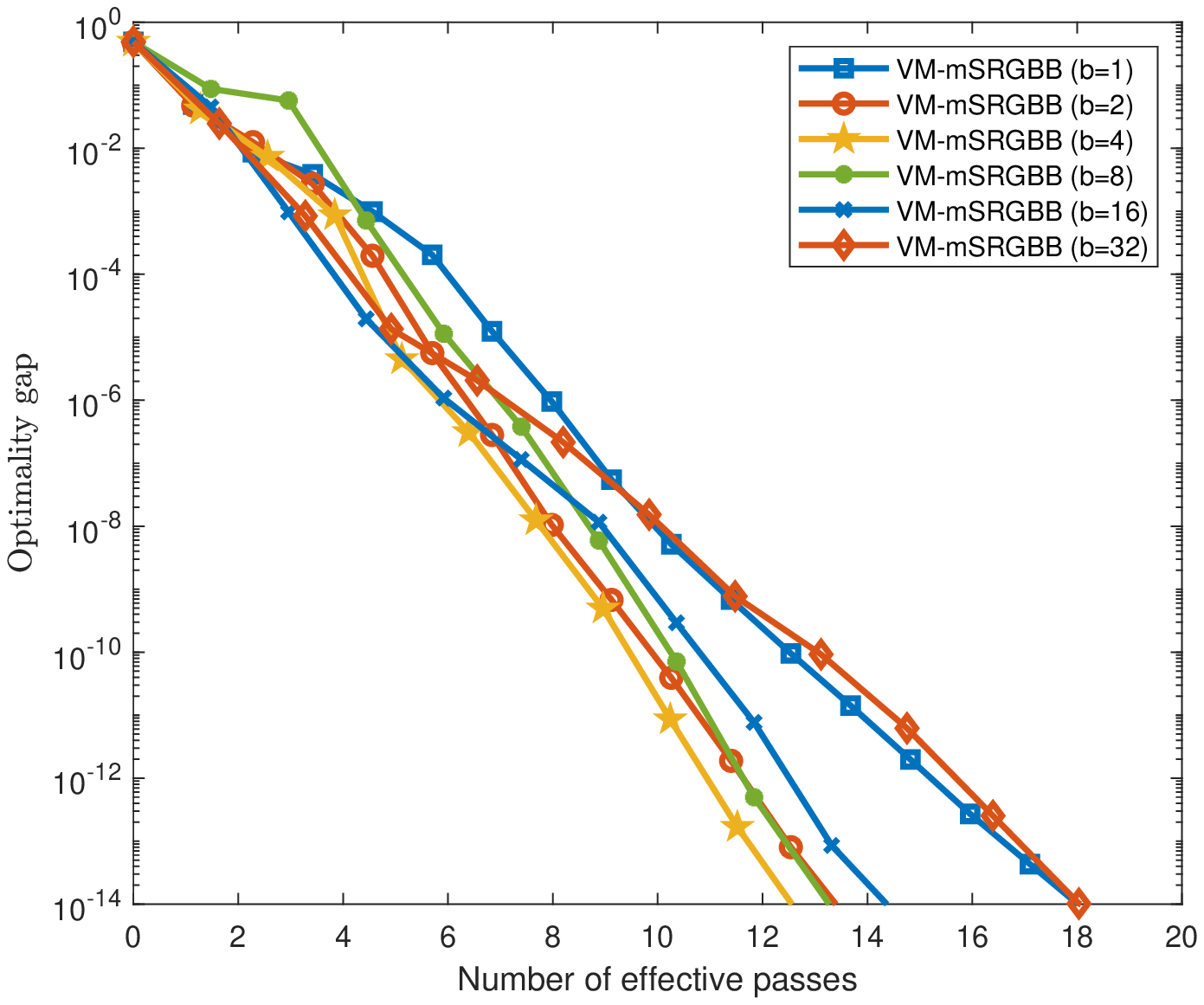}}
  \hspace{3mm}
  \subfigure[rcv1]{
    \label{fig:subfig:rcv1b} %% label for first subfigure
    \includegraphics[width=1.6in,height=1.3in]{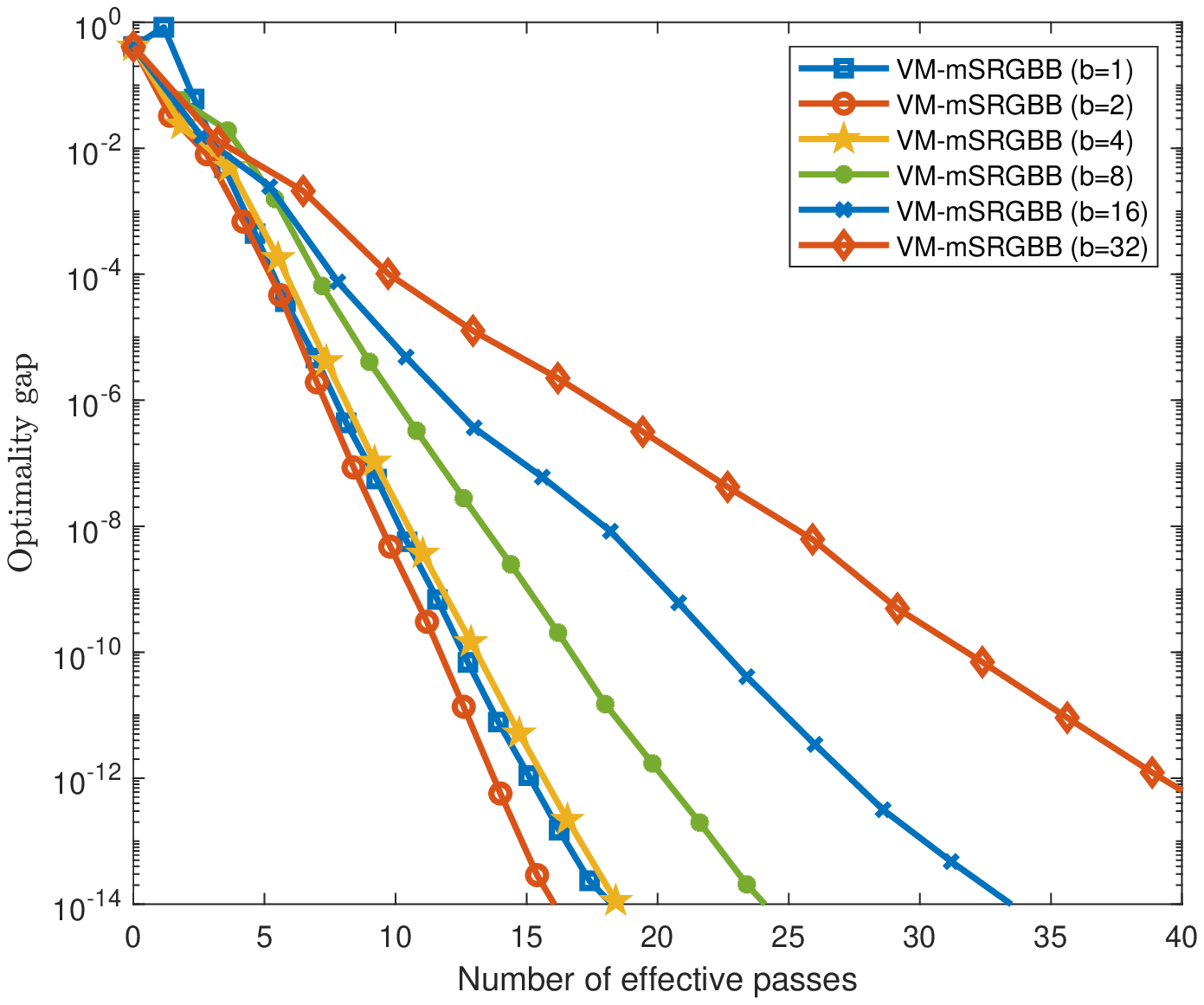}}
  \subfigure[real-sim]{
    \label{fig:subfig:real-simc} %% label for first subfigure
    \includegraphics[width=1.6in,height=1.3in]{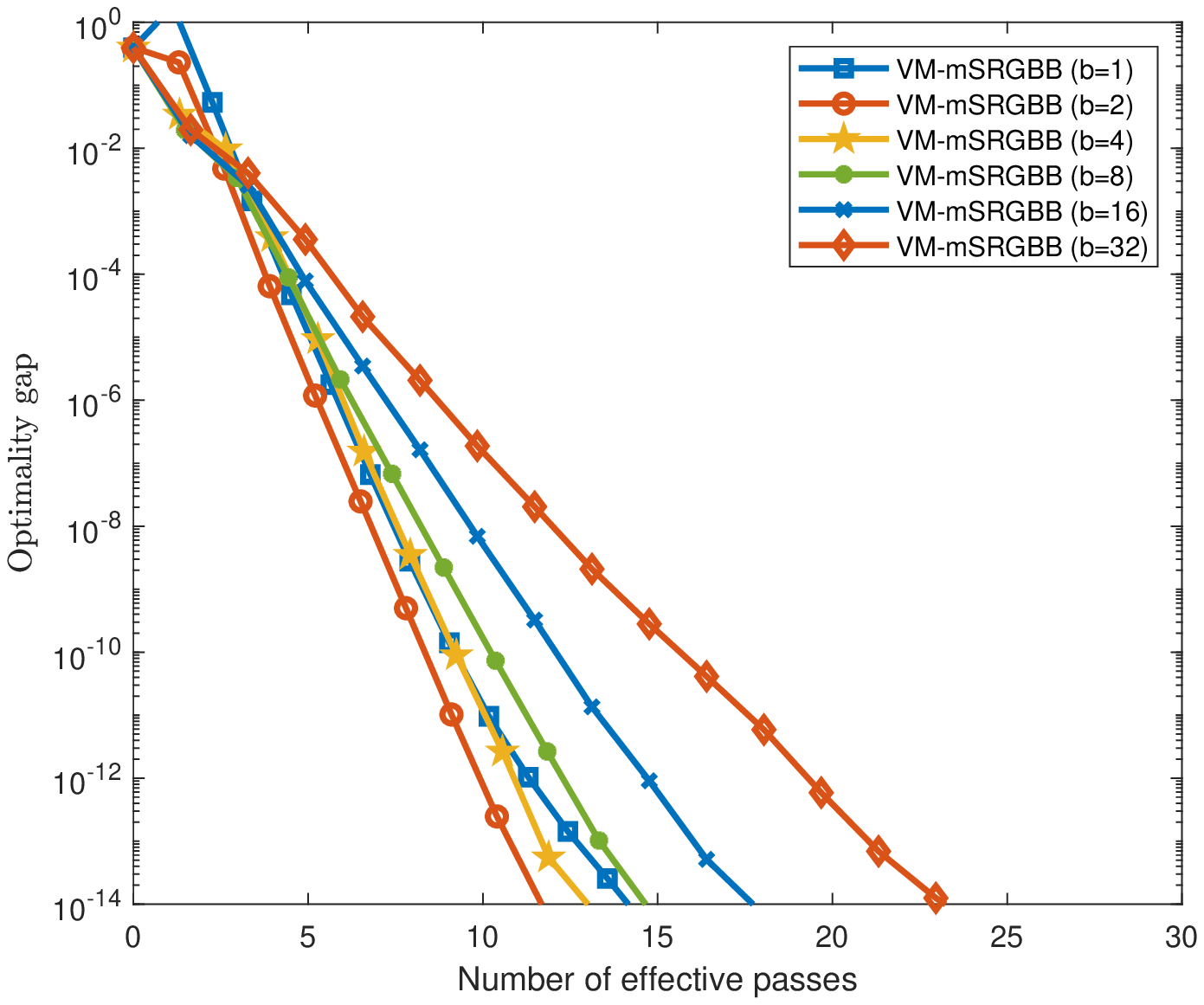}}
  \hspace{3mm}
  \subfigure[covtype]{
    \label{fig:subfig:covtyped} %% label for first subfigure
    \includegraphics[width=1.6in,height=1.3in]{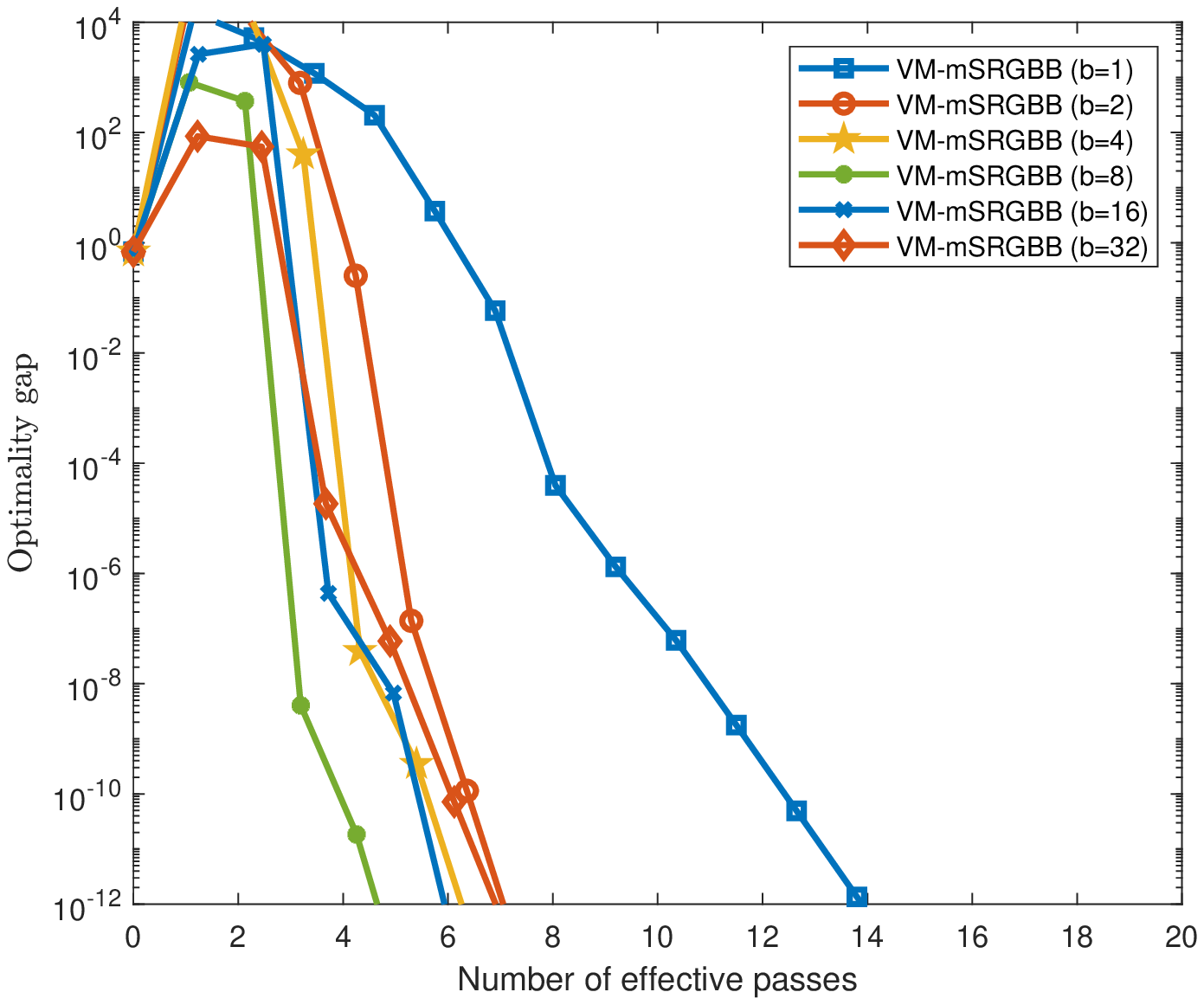}}
  \caption{Comparison of VM-mSRGBB with different mini-batch sizes. }\label{VSmini}
\end{figure}

Fig. \ref{VSmini} illustrates the results of VM-mSRGBB under various mini-batch sizes $b$ on the four data sets. We can see that compared with $b=1$, VM-mSRGBB has better or comparable performance by increasing the mini-batch size to $b = 2, 4, 8, 16$.

\subsection{Comparison with other algorithms}\label{ne3}
In this part, we conduct experiments on VM-mSRGBB in comparison with four modern mini-batch proximal stochastic gradient methods, which are specified as follows:
\begin{description}
  \item[(1)] mS2GD: mS2GD is a mini-batch proximal version of S2GD \cite{konevcny2013semi} to deal with nonsmooth problems. In mS2GD, a constant stepsize was used.
  \item[(2)] mS2GD-BB: mS2GD-BB uses the BB method to compute stepsizes for mS2GD.
%  mS2GD-BB was run with the best-tuned parameters.
  \item[(3)] mSARAH: mSARAH is a mini-batch proximal variant of   stochastic recursive gradient algorithm proposed in \cite{nguyen2017sarah}. In mSARAH, a constant stepsize was used.
  \item[(4)] mSARAH-BB: mSARAH-BB is a mini-batch variant of SARAH-BB \cite{liu2020a}.
%  We used the mini-batch variant of BB stepsize.
\end{description}

\begin{figure}
%\centering
%\begin{minipage}[!t]{0.5\linewidth}
  \subfigure[ijcnn1]{
    \label{fig:subfig:ijcnn1a} %% label for first subfigure
    \includegraphics[width=1.6in,height=1.3in]{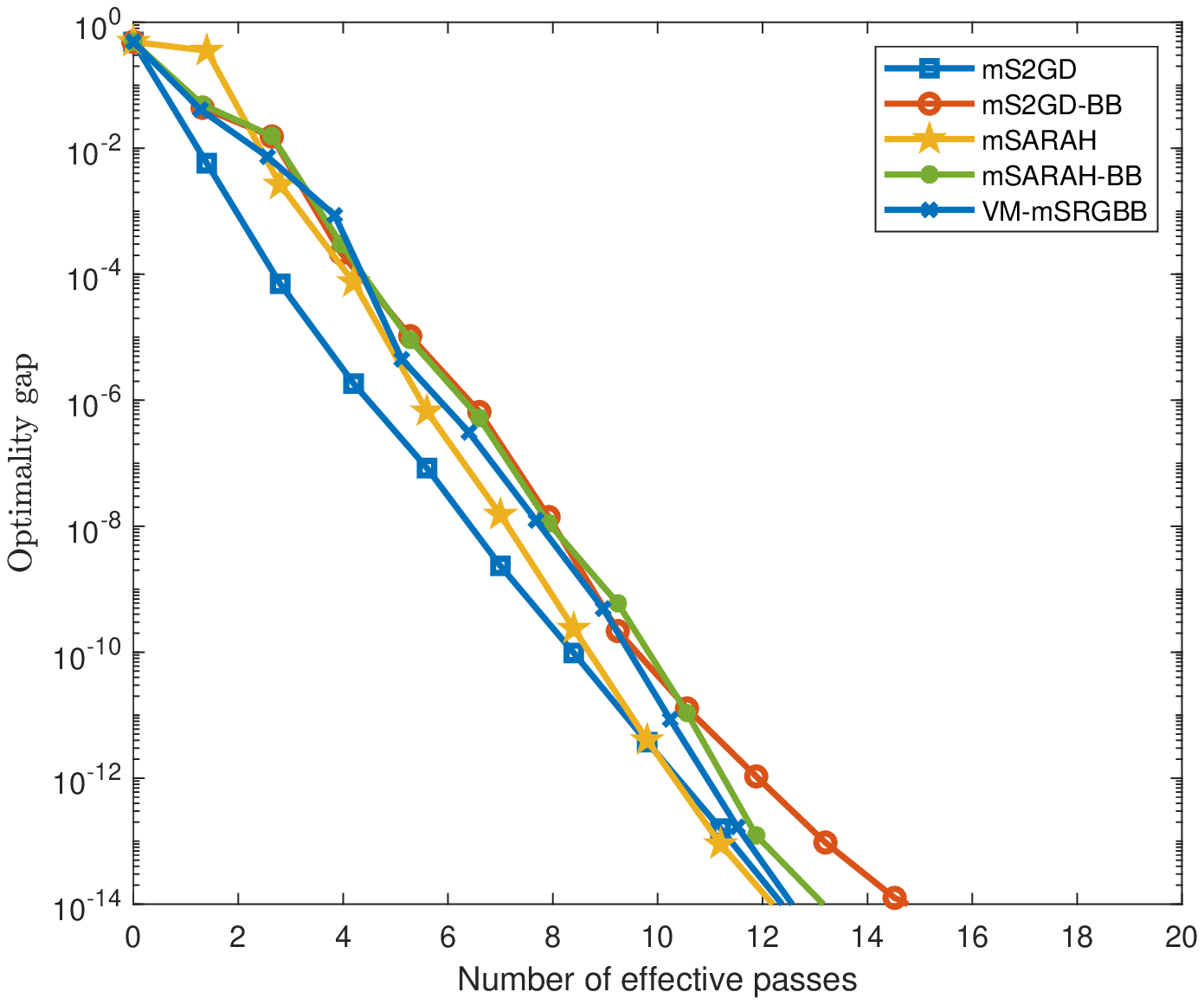}}
  \hspace{3mm}
  \subfigure[rcv1]{
    \label{fig:subfig:rcv1b} %% label for first subfigure
    \includegraphics[width=1.6in,height=1.3in]{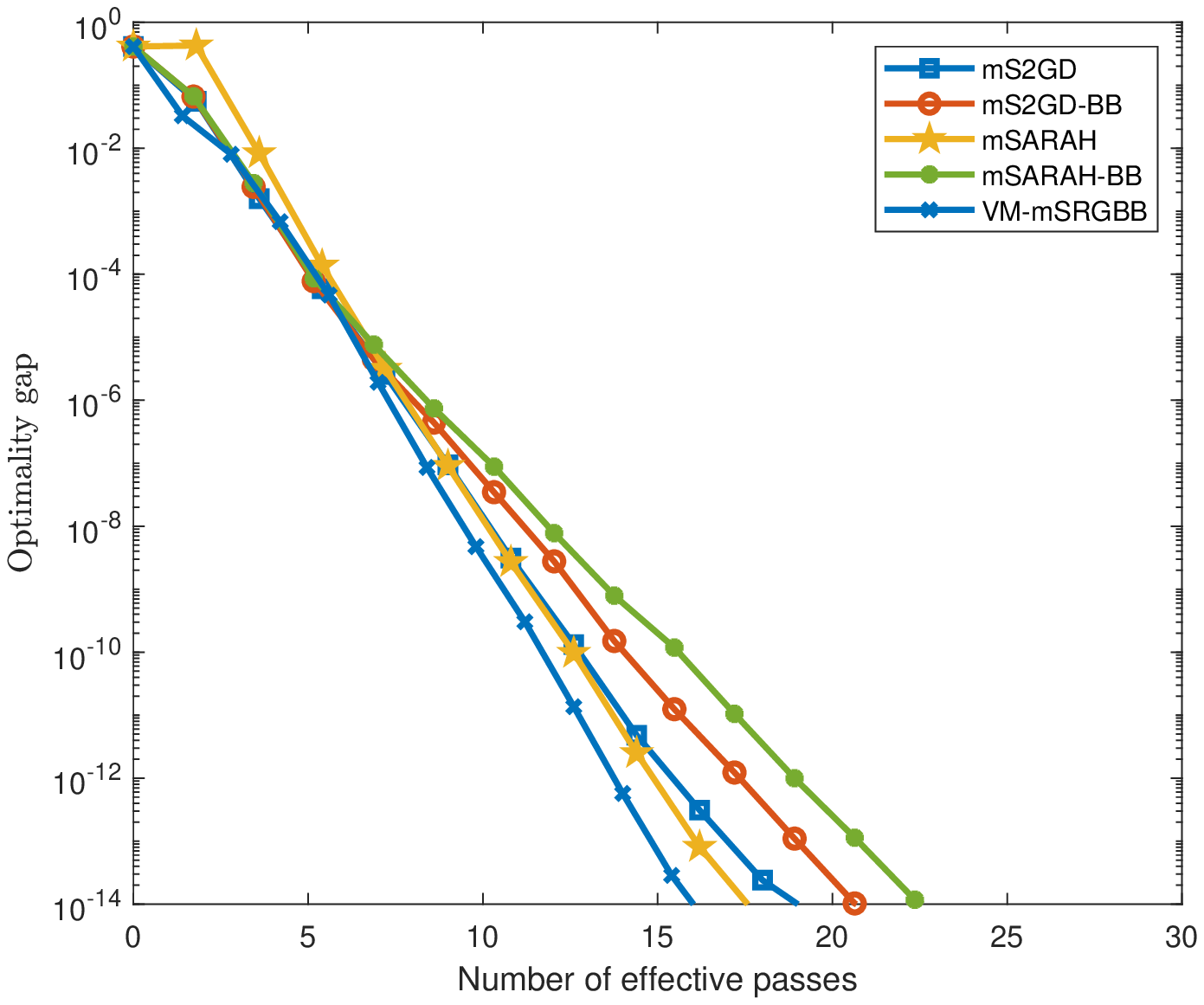}}
  \subfigure[real-sim]{
    \label{fig:subfig:real-simc} %% label for first subfigure
    \includegraphics[width=1.6in,height=1.3in]{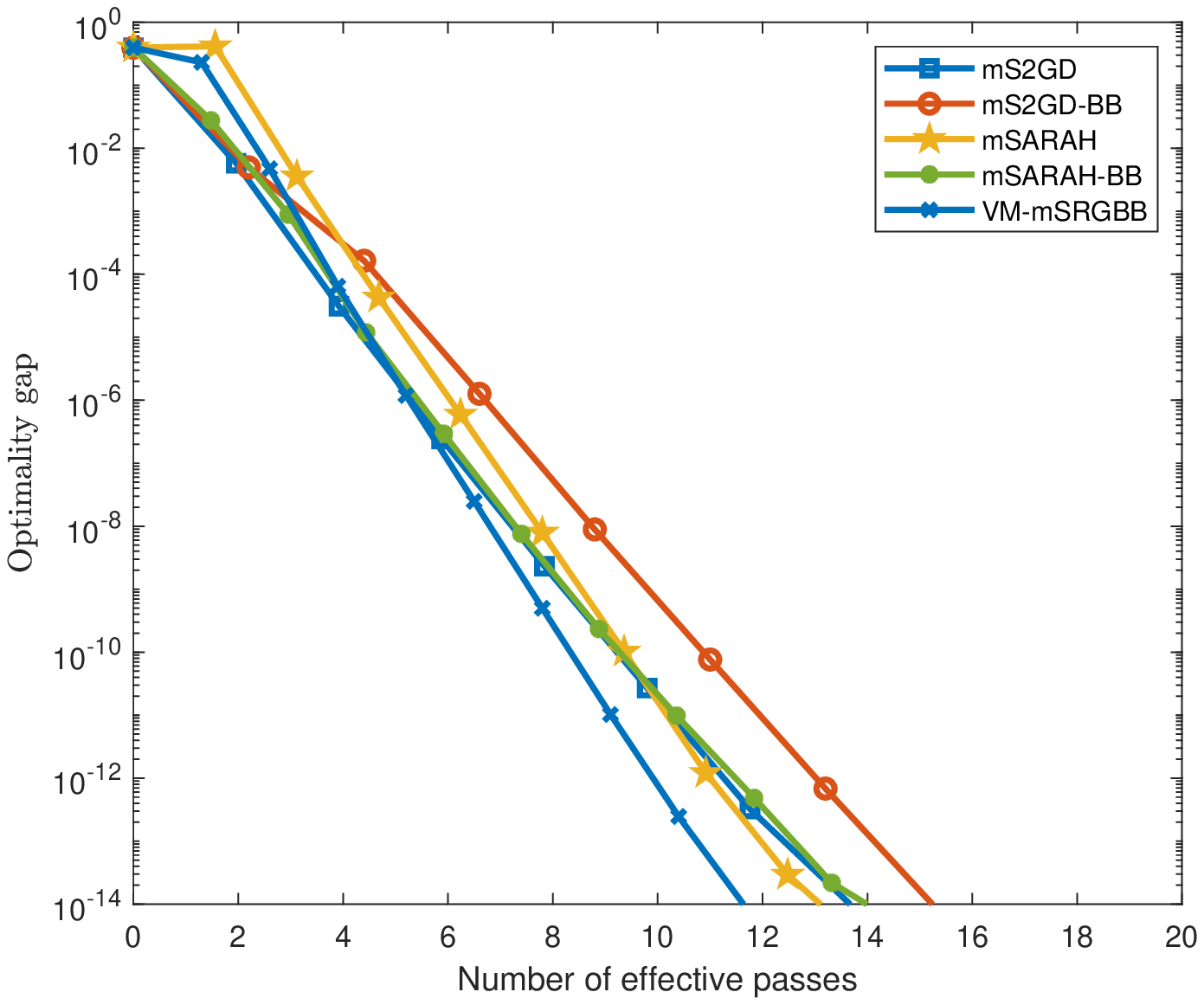}}
    \hspace{3mm}
  \subfigure[covtype]{
    \label{fig:subfig:covtyped} %% label for first subfigure
    \includegraphics[width=1.6in,height=1.3in]{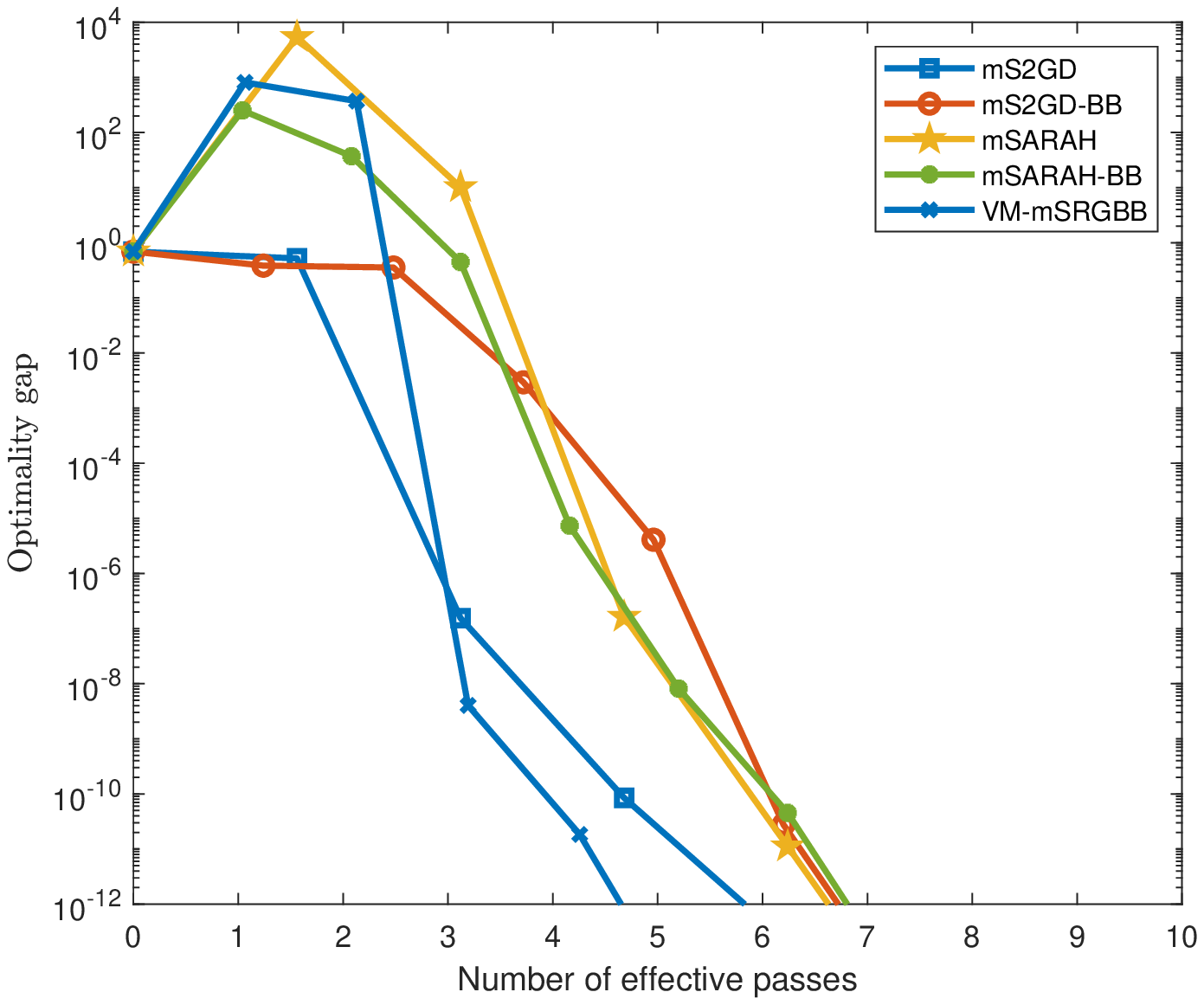}}
  \caption{Comparison of VM-mSRGBB and other modern methods.}\label{VSothers}
\end{figure}

For the above four methods, we used $b=8$. The choices of parameters employed by VM-mSRGBB are given in Table \ref{tab5}. Fig. \ref{VSothers} demonstrates that our VM-mSRGBB is superior to the compared algorithms on the four data sets.
\begin{table}[!t]
\renewcommand{\arraystretch}{1.3}
\caption{Best choices of parameters in VM-mSRGBB}
\label{tab5}
\centering
\begin{tabular}{|c|c|c|c|c|}
  \hline
  Parameter & ijcnn1 & rcv1 & real-sim & covtype \\
  \hline
  $b$ & 4 & 2 & 2 & 8 \\
  \hline
  $m$ & 0.07$n$ & 0.2$n$ & 0.15$n$ & 0.008$n$\\
  \hline
\end{tabular}
\end{table}

\section{Conclusion}\label{conclusion}
Based on a newly derived diagonal BB stepsizes for updating the metric, we proposed a proximal stochastic recursive gradient method named VM-mSRGBB to minimize the composition of two convex functions.
Linear convergence of VM-mSRGBB was established under mild conditions for strongly convex, non-strongly convex and convex cases, respectively. Numerical comparisons of VM-mSRGBB and recent successful stochastic variance reduced gradient methods and mini-batch proximal stochastic methods on some real data sets highly suggest the potential benefits of our VM-mSRGBB method for composition optimization problems arising in machine learning.

\appendices

\section{Proof of Lemma \ref{vmlem6}}\label{vmapp4}
We only need to consider the nontrivial case $w\neq w'$. Denoting $\mathbf{u} = \textrm{prox}_{R}^{A^{-1}}(w)$ and $\mathbf{v} = \textrm{prox}_R^{A^{-1}}(w')$. It follows from Lemma \ref{vmlem2} that
\begin{equation*}\label{eqvms2.3}
  A^{-1}(w-\mathbf{u}) \in \partial R(\mathbf{u}), \quad A^{-1}(w'-\mathbf{v}) \in \partial R(\mathbf{v}).
\end{equation*}
By the definition of subdifferential, we have
\begin{equation*}\label{eqvms2.4}
  R(\mathbf{u}) \geq R(\mathbf{v}) + (A^{-1}(w'-\mathbf{v}))^T(\mathbf{u} - \mathbf{v}),
\end{equation*}
\begin{equation*}\label{eqvms2.5}
  R(\mathbf{v}) \geq R (\mathbf{u}) + (A^{-1}(w-\mathbf{u}))^T(\mathbf{v} - \mathbf{u}).
\end{equation*}
Summing the above two inequalities to get
\begin{align*}\label{eqvms2.6}
  0 &\geq \Big(A^{-1}\big((w'-\mathbf{v}) - (w-\mathbf{u})\big)\Big)^T\Big(\mathbf{u}-\mathbf{v}\Big) \nonumber\\
   &=\Big(A^{-1}\big((w'-w)+(\mathbf{u}-\mathbf{v})\big)\Big)  ^T\Big(\mathbf{u}-\mathbf{v}\Big),
\end{align*}
which results in,
\begin{align*}
  \|\mathbf{u}-\mathbf{v}\|_{A^{-1}}^2 &\leq (A^{-1}(w-w'))^T(\mathbf{u}-\mathbf{v})\\
   %&=&  \|A^{-1/2}(\mathbf{u}-\mathbf{v})\|_2^2 \\
   &= (A^{-1/2}(w-w'))^T (A^{-1/2}(\mathbf{u}-\mathbf{v}))\\
   &\leq \|A^{-1/2}(w-w')\|_2 \cdot \|A^{-1/2}(\mathbf{u}-\mathbf{v})\|_2,
\end{align*}
where the first equality holds due to the symmetry and positive definiteness of $A$ while the last inequality follows from the Cauchy-Schwarz inequality. By squaring the above inequality, we obtain
\begin{align*}
  &\|\mathbf{u}-\mathbf{v}\|_{A^{-1}}^2\cdot\|\mathbf{u}-\mathbf{v}\|_{A^{-1}}^2\\
  \leq& \|A^{-1/2}(w-w')\|_2^2 \cdot \|A^{-1/2}(\mathbf{u}-\mathbf{v})\|_2^2 \\
   =& \|w-w'\|_{A^{-1}}^2 \cdot \|\mathbf{u}-\mathbf{v}\|_{A^{-1}}^2.
\end{align*}
Since $w\neq w'$, we know that $\|\mathbf{u}-\mathbf{v}\|_{A^{-1}}^2\neq0$. We complete the proof by dividing both sides of the above inequality by $\|\mathbf{u}-\mathbf{v}\|_{A^{-1}}^2$.

\section{Proof of Lemma \ref{vmlem5}}\label{vmapp1}
Since
$$w_{t+1}^k = \arg\min_y\big\{R(y) + \frac{1}{2}\|y - (w_t^k - U_kv_t^k)\|_{U_k^{-1}}^2\big\},$$
by Lemma \ref{vmlem2}, we get
\begin{equation*}
U_k^{-1}\big((w_t^k - U_kv_t^k)-w_{t+1}^k\big) \in \partial R(w_{t+1}^k),
\end{equation*}
which implies that there exists $\varphi\in \partial R(w_{t+1}^k)$ such that
\begin{equation*}
  U_k^{-1}\big(w_{t+1}^k - (w_t^k - U_kv_t^k)\big) + \varphi =0.
\end{equation*}
This together with \eqref{eqvmspm2} gives
\begin{equation*}\label{eqvmc5}
  v_t^k+\varphi = g_t^k.
\end{equation*}
Then
\begin{equation}\label{eqvmc1}
  (w_*-w_{t+1}^k)^T(v_t^k+\varphi) = (w_*-w_{t+1}^k)^Tg_t^k.
\end{equation}
From the convexity of $F(w)$ and $R(w)$, we get
\begin{align}\label{eqvmc6}
  P(w_*) \geq& F(w_t^k) + \nabla F(w_t^k)^T(w_* - w_t^k) \nonumber\\
  & + R(w_{t+1}^k) + \varphi^T(w_* - w_{t+1}^k).
\end{align}
It follows from the Lipschitz continuity of $\nabla F(w)$ that
\begin{align}\label{eqvmc7}
  &F(w_t^k) \nonumber\\
  \geq& F(w_{t+1}^k) - \nabla F(w_t^k)^T(w_{t+1}^k - w_t^k)
  - \frac{L}{2}\|w_{t+1}^k - w_t^k\|_2^2 \nonumber\\
  \geq& F(w_{t+1}^k) - \nabla F(w_t^k)^T(w_{t+1}^k - w_t^k) - \frac{L_{\Omega}}{2}\|w_{t+1}^k - w_t^k\|_2^2,
\end{align}
where the second inequality is due to the fact $0< L\leq L_{\Omega}$. Combining \eqref{eqvmc6} and \eqref{eqvmc7}, we have
\begin{align}\label{eqvmc3}
  & P(w_*) \nonumber\\
  \geq & F(w_{t+1}^k) - \nabla F(w_t^k)^T(w_{t+1}^k - w_t^k) + \nabla F(w_t^k)^T(w_* - w_t^k)\nonumber\\
  & +  R(w_{t+1}^k) + \varphi^T(w_* - w_{t+1}^k) - \frac{L_{\Omega}}{2}\|w_{t+1}^k - w_t^k\|_2^2\nonumber\\
  =& P(w_{t+1}^k)  + \nabla F(w_t^k)^T(w_* - w_{t+1}^k)+ \varphi^T(w_* - w_{t+1}^k)\nonumber\\
  &- \frac{L_{\Omega}}{2}\|w_{t+1}^k - w_t^k\|_2^2\nonumber\\
   \geq& P(w_{t+1}^k)  + \nabla F(w_t^k)^T(w_* - w_{t+1}^k) + \varphi^T(w_* - w_{t+1}^k)\nonumber\\
   & -\frac{1}{2}\|g_t^k\|_{U_k}^2,
\end{align}
where the first equality follows from the definition of $P(w)$ and the last inequality holds by \eqref{eqvmspm2} and $0\prec U_k \preceq 1/L_{\Omega}I$. Collecting all inner products on the right-hand side of \eqref{eqvmc3}, we obtain
\begin{align}\label{eqvmc2}
 & \nabla F(w_t^k)^T(w_* - w_{t+1}^k) + \varphi^T(w_* - w_{t+1}^k)  \nonumber\\
   =&  (w_* - w_{t+1}^k)^T(\delta_t^k+v_t^k) +(w_* - w_{t+1}^k)^T\varphi\nonumber\\
   =& (w_* - w_{t+1}^k)^T\delta_t^k + (w_* - w_{t+1}^k)^T(v_t^k + \varphi)\nonumber\\
   = & (w_* - w_{t+1}^k)^T\delta_t^k + (w_* - w_{t+1}^k)^Tg_t^k \nonumber\\
   = & (w_* - w_{t+1}^k)^T\delta_t^k + (w_* - w_t^k + w_t^k-w_{t+1}^k)^Tg_t^k \nonumber\\
   =& (w_* - w_{t+1}^k)^T\delta_t^k + (w_* - w_t^k)^Tg_t^k + (g_t^k)^TU_kg_t^k \nonumber\\
   =& (w_* - w_{t+1}^k)^T\delta_t^k + (w_* - w_t^k)^Tg_t^k + \|g_t^k\|_{U_k}^2,
\end{align}
where the first equality follows from the definition of $\delta_t^k$, and the third and fifth equalities are derived from \eqref{eqvmc1} and \eqref{eqvmspm2}, respectively. Applying \eqref{eqvmc2} to \eqref{eqvmc3}, we get
\begin{align*}
  P(w_*) \geq&  P(w_{t+1}^k) + \frac{1}{2}\|g_t^k\|_{U_k}^2 + (w_* - w_{t+1}^k)^T\delta_t^k \\
  & + (w_* - w_t^k)^Tg_t^k.
\end{align*}
Then the desired result is obtained.

\section{Proof of Lemma \ref{vmlem10.1}}\label{vmapp2}
Consider $v_t^k$ defined in \eqref{eqvmvt1.1}. Conditioned on $\mathcal{F}_t = \sigma(w_0^k,i_1,\ldots,i_{t-1})$, we take expectation with respect to $i_t$ and obtain
\begin{equation}\label{eqvms1}
  \mathbb{E}\Big[\frac{\nabla f_{i_t}(w_t^k)}{nq_{i_t}} |\mathcal{F}_t\Big] = \sum_{i=1}^n\frac{q_i}{nq_i}\nabla f_i(w_t^k) = \nabla F(w_t^k).
\end{equation}
Similarly we have
\begin{equation}\label{eqvms2}
  \mathbb{E}\Big[\frac{\nabla f_{i_t}(w_{t-1}^k)}{nq_{i_t}} |\mathcal{F}_t\Big] = \nabla F(w_{t-1}^k).
\end{equation}
Then we obtain
\begin{align*}
  & \mathbb{E}\Big[\|v_t^k-\nabla F(w_t^k)\|_2^2|\mathcal{F}_t\Big] \\
 = &  \mathbb{E}\Big[\|\frac{\nabla f_{i_t}(w_t^k) - \nabla f_{i_t}(w_{t-1}^k)}{nq_{i_t}} - \big(\nabla F(w_t^k) - \nabla F(w_{t-1}^k)\big) \\
 & + \big(v_{t-1}^k - \nabla F(w_{t-1}^k)\big)\|_2^2|\mathcal{F}_t\Big] \\
 =& \mathbb{E}\Big[\|\frac{\nabla f_{i_t}(w_t^k) - \nabla f_{i_t}(w_{t-1}^k)}{nq_{i_t}}\|_2^2|\mathcal{F}_t\Big]\\
 & - \|\nabla F(w_t^k) - \nabla F(w_{t-1}^k)\|_2^2 + \|v_{t-1}^k - \nabla F(w_{t-1}^k)\|_2^2\\
 =&  \mathbb{E}\Big[\|\frac{\nabla f_{i_t}(w_t^k) - \nabla f_{i_t}(w_{t-1}^k)}{nq_{i_t}}\|_2^2|\mathcal{F}_t\Big]\\
 & - 2\big(\nabla F(w_t^k) -v_{t-1}^k\big)^T\big(v_{t-1}^k - \nabla F(w_{t-1}^k)\big) \\
 & - \|\nabla F(w_t^k) - v_{t-1}^k\|_2^2,
\end{align*}
where the second equality follows from \eqref{eqvms1} and \eqref{eqvms2}.

By taking total expectation over the entire history in the $k$-th outer loop, we obtain
\begin{align*}
   & \mathbb{E}\big[\|v_t^k-\nabla F(w_t^k)\|_2^2\big]  \\
   =&  \mathbb{E}\big[\mathbb{E}[\|\nabla F(w_t^k) - v_t^k\|_2^2|\mathcal{F}_t]\big]\\
   = & \mathbb{E}\big[\|\frac{\nabla f_{i_t}(w_t^k) - \nabla f_{i_t}(w_{t-1}^k)}{nq_{i_t}}\|_2^2\big] - \mathbb{E}[\|\nabla F(w_t^k) - v_{t-1}^k\|_2^2]\\
   \leq&  \mathbb{E}\big[\|\frac{\nabla f_{i_t}(w_t^k) - \nabla f_{i_t}(w_{t-1}^k)}{nq_{i_t}}\|_2^2\big]\\
   \leq& \mathbb{E}\big[\frac{L_{i_t}^2}{n^2q_{i_t}^2}\|w_t^k - w_{t-1}^k\|_2^2\big] \\
   \leq& L_{\Omega}^2 \mathbb{E}\big[ \|w_t^k - w_{t-1}^k\|_2^2\big],
\end{align*}
where the second equality holds due to \eqref{eqvmtotal2}, the second inequality follows from the Lipschitz continuity of $\nabla f_i$, and the last inequality is due to the fact that $L_{\Omega} \geq \frac{L_i}{nq_i}$ for $i =1,2,\ldots,n$.

\section{Proof of Lemma \ref{vmlem10}}\label{vmapp3}
We definite $G_i = \frac{\nabla f_{i}(w_t^k) - \nabla f_{i}(w_{t-1}^k)}{nq_{i}} + v_{t-1}^k$, then $v_t^k$ in \eqref{eqvmvt1} can be written as
\begin{equation*}
  v_t^k = \frac{1}{b}\sum_{i\in I_t}\Big(\frac{\nabla f_{i}(w_t^k) - \nabla f_{i}(w_{t-1}^k)}{nq_{i}} + v_{t-1}^k\Big) = \frac{1}{b}\sum_{i\in I_t} G_i.
\end{equation*}
Conditioned on $\mathcal{F}_t = \sigma(w_0^k,I_1,\ldots,I_{t-1})$, we take expectation with respect to $I_t$ and get
\begin{align*}
  & \mathbb{E}\big[\|v_t^k - \nabla F(w_t^k)\|_2^2|\mathcal{F}_t\big] \\
  =& \frac{1}{b^2} \mathbb{E}\big[\|\sum_{i\in I_t}(G_i- \nabla F(w_t^k))\|_2^2|\mathcal{F}_t\big] \\
  =& \frac{1}{b^2}\mathbb{E}\big[\|\sum_{i\in S_1}(G_i- \nabla F(w_t^k))+ (G_{I_t/S_1}-\nabla F(w_t^k))\|_2^2|\mathcal{F}_t\big] \\
  =&  \frac{1}{b^2}\mathbb{E}\big[\|\sum_{i\in S_1}(G_i- \nabla F(w_t^k))\|_2^2|\mathcal{F}_t\big] \\
  &+ \frac{1}{b^2}\mathbb{E}\big[\|G_{I_t/S_1}-\nabla F(w_t^k)\|_2^2|\mathcal{F}_t\big]\\
  & + \frac{2}{b^2}\mathbb{E}\big[\big(\sum_{i\in S_1}(G_i- \nabla F(w_t^k))\big)^T\big(G_{I_t/S_1}-\nabla F(w_t^k)\big)|\mathcal{F}_t\big],
\end{align*}
where $S_1\subset I_t$ and the number of elements in the set $I_t/S_1$ is 1. By taking expectation over the entire history and applying the above inequality recursively, we obtain
\begin{align*}
  &\mathbb{E}[\|v_t^k - \nabla F(w_t^k)\|_2^2] \\
  =& \frac{1}{b^2}\mathbb{E}\big[\|\sum_{i\in S_1}(G_i- \nabla F(w_t^k))\|_2^2\big] \\
  &+ \frac{1}{b^2}\mathbb{E}\big[\|G_{I_t/S_1}-\nabla F(w_t^k)\|_2^2\big]\\
  =& \frac{1}{b^2}\sum_{i\in I_t}\mathbb{E}\big[\|G_i- \nabla F(w_t^k)\|_2^2\big]\\
  \leq& \frac{L_{\Omega}^2}{b} \mathbb{E}\big[ \|w_t^k - w_{t-1}^k\|_2^2\big],
\end{align*}
where the first equality holds due to the fact $\mathbb{E}[G_i] = \mathbb{E}[\nabla F(w_t^k)]$, which follows from \eqref{eqvmtotal2} with $b=1$. In the last inequality we use Lemma \ref{vmlem10.1}.

\section*{Acknowledgment}
This research is supported by the Chinese NSF grants (nos. 11671116, 11701137, 12071108, 11631013, 11991020 and 12021001), the Major Research Plan of the NSFC (no. 91630202), and Beijing Academy of Artificial Intelligence (BAAI).

%\section*{References}
%

%\bibliography{IEEEabrv,VM-mSRGBB}

%\appendix[IEEEproof of the Zonklar Equations]
% or


%\appendices
%\section{IEEEproof of the First Zonklar Equation}
%% you can choose not to have a title for an appendix
%% if you want by leaving the argument blank
%\section{}

%\section*{Acknowledgment}
%This work is supported by the Chinese NSF grants (nos. 11671116 and 11271107) and the Major Research Plan of the NSFC (no. 91630202).
%
%
%The authors would like to thank anonymous reviewers and the handling editor for their constructive comments and suggestions.

\ifCLASSOPTIONcaptionsoff
  \newpage
\fi



\bibliographystyle{IEEEtran}

\end{document}